\documentclass[oneeqnum,onetabnum,onefignum,onealgnum,onethmnum]{siamart251216}

\newcommand{\TheTitle}{%
  Multigrid with Linear Storage Complexity
}

\newcommand{\TheFunding}{%
  The author/s would like to thank the NHR-Verein e.V. (\url{www.nhr-verein.de}) for supporting this work/project within the NHR Graduate School of National High Performance Computing (NHR).
  The authors gratefully acknowledge support through the joint project CoMPS (\url{https://gauss-allianz.de/en/project/title/CoMPS}) (grant \texttt{16ME0647K}) funded by the German Federal Ministry of Research, Technology and Space.
  The authors gratefully acknowledge the Gauss Centre for Supercomputing e.V. (\url{https://www.gauss-centre.eu}) for funding this project by providing computing time on the GCS Supercomputer SuperMUC-NG and friendly user access during the pilot operation of the GCS Supercomputer SuperMUC-NG Phase 2 at Leibniz Supercomputing Centre (\url{https://www.lrz.de}).
}

\author{%
  Daniel Bauer\,\orcidlink{0009-0008-9397-2883}\thanks{Friedrich-Alexander-Universität Erlangen-Nürnberg (FAU), Erlangen National High Performance Computing Center (NHR@FAU), (\email{daniel.j.bauer@fau.de}).}\and%
  Nils Kohl\,\orcidlink{0000-0003-4797-0664}\thanks{Dept.\ of Earth and Environmental Sciences, Ludwig-Maximilians-Universität München (LMU).}\and%
  Stephen F. McCormick\,\orcidlink{0000-0002-6327-790X}\thanks{University of Colorado at Boulder, Boulder, CO.}\and%
  Rasmus~Tamstorf\,\orcidlink{0000-0001-7649-8587}\thanks{Independent researcher.}%
}
\title{{\TheTitle}\thanks{\TheFunding}}

\usepackage{amsmath, amssymb}
\usepackage{siunitx}
\usepackage{pseudo}
\pseudoset{
  ct-left={$\triangleright$ },
  ct-right={},
  indent-mark,
  indent-mark-color=black,
  indent-mark-width=.4pt,
  indent-text={\pseudofont\kw{if}\ },
}
\pseudodefinestyle{fullwidth}{
  begin-tabular = \tabularx{0.98\linewidth}[t]{>{\pseudohpad}\pseudolabelalign>{\leavevmode\pseudosetup}X<{\pseudohpad}},
  end-tabular = \endtabularx
}
\pseudodefinestyle{commentcol}{
  begin-tabular = \tabularx{0.98\linewidth}[t]{>{\pseudohpad}\pseudolabelalign>{\leavevmode\pseudosetup}X>{\leavevmode\pseudosetup}p{0.35\linewidth}<{\pseudohpad}},
  end-tabular = \endtabularx
}
\makeatletter
\NewDocumentCommand \pseudoanchor { m } {%
  \tikz[baseline, overlay, remember picture]
    \node[anchor=base, inner sep=0] (#1) {\@arstrut};%
  \ignorespaces
}
\makeatother

\usepackage{pgfplots, tikz}
\usetikzlibrary{decorations.pathreplacing,matrix}
\usepgfplotslibrary{groupplots}
\tikzset{every picture/.append style={remember picture}}
\pgfplotsset{compat=1.18}

\pgfplotsset{
  colormap name = viridis,
  grid = both,
  grid style = {line width=.2pt, draw=gray!20},
  major grid style = {line width=.2pt,draw=gray!50},
  legend style = {
    font = \footnotesize,
    draw = none,
    inner sep = 0pt,
  },
  label style = {
    font = \footnotesize,
    align = center,
  },
  ticklabel style = {
    font = \footnotesize,
  },
  title style = {
    font = \footnotesize,
    yshift = 0pt,
  },
}

\definecolor{pm-indg}{HTML}{332288}
\definecolor{pm-cyan}{HTML}{88ccee}
\definecolor{pm-teal}{HTML}{44aa99}
\definecolor{pm-gren}{HTML}{117733}
\definecolor{pm-oliv}{HTML}{999933}
\definecolor{pm-sand}{HTML}{ddcc77}
\definecolor{pm-rose}{HTML}{cc6677}
\definecolor{pm-wine}{HTML}{882255}
\definecolor{pm-prpl}{HTML}{aa4499}

\definecolor{pm-rb09}{HTML}{882e72}
\definecolor{pm-rb10}{HTML}{1965b0}
\definecolor{pm-rb15}{HTML}{4eb265}
\definecolor{pm-rb24}{HTML}{e8601c}
\definecolor{pm-rb26}{HTML}{dc050c}
\definecolor{pm-rb27}{HTML}{a5170e}
\definecolor{pm-rb28}{HTML}{72190e}

\definecolor{db-redd}{HTML}{700000}
\definecolor{db-prpl}{HTML}{833399}

\pgfplotscreateplotcyclelist{mymarklist}{
  {mark=*          ,mark options={scale=0.5}},
  {mark=square*    ,mark options={scale=0.5}},
  {mark=triangle*  ,mark options={scale=0.5}},
  {mark=halfcircle*,mark options={scale=0.5}},
  {mark=halfsquare*,mark options={scale=0.5}},
  {mark=o          ,mark options={scale=0.5}},
  {mark=square     ,mark options={scale=0.5}},
  {mark=triangle   ,mark options={scale=0.5}},
  {mark=pentagon   ,mark options={scale=0.5}},
  {mark=diamond    ,mark options={scale=0.5}}
}

\pgfplotsset{
  cycle list = {
    {pm-rb10,mark=*          ,mark options={scale=0.5}},
    {pm-rb24,mark=square*    ,mark options={scale=0.5}},
    {pm-gren,mark=triangle*  ,mark options={scale=0.5}},
    {pm-rb09,mark=halfcircle*,mark options={scale=0.5}},
    {black  ,mark=halfsquare*,mark options={scale=0.5}},
    {pm-rb10,mark=o          ,mark options={scale=0.5}},
    {pm-rb24,mark=square     ,mark options={scale=0.5}},
    {pm-gren,mark=triangle   ,mark options={scale=0.5}},
    {pm-rb09,mark=pentagon   ,mark options={scale=0.5}},
    {black  ,mark=diamond    ,mark options={scale=0.5}}
  },
}

\usepackage[caption=false]{subfig}
\usepackage{booktabs,multirow,tabularx}

\usepackage{orcidlink}
\crefname{pseudoline}{line}{lines}
\newsubfloat[position=top]{algorithm}
\crefname{subalgorithm}{Algorithm}{Algorithms}
\newsiamremark{remark}{Remark}

\usepackage{glossaries, glossaries-prefix}
\newacronym{bfp}{BFP}{block floating point}
\newacronym{bfs}{BFS}{breadth-first search}
\newacronym{bpx}{BPX}{Bramble-Pasciak-Xu}
\newacronym{cg}{CG}{conjugate gradient}
\newacronym{cs}{CS}{correction scheme}
\newacronym{cmg}{CMG}{compact multigrid}
\newacronym{dfs}{DFS}{depth-first search}
\newacronym[longplural={degrees of freedom}]{dof}{DoF}{degree of freedom}
\newacronym{fe}{FE}{finite element}
\newacronym{fas}{FAS}{full approximation scheme}
\newacronym{fem}{FEM}{finite element method}
\newacronym{flop}{flop}{floating point operation}
\newacronym[prefix = {an\space}, prefixfirst = {a~}]{fmg}{FMG}{full multigrid}
\newacronym{fp}{FP}{floating point}
\newacronym{gemv}{gemv}{general matrix-vector product}
\newacronym{ir}{IR}{iterative refinement}
\newacronym{lsb}{LSB}{least-significant bit}
\newacronym[prefix = {an\space}, prefixfirst = {a~}]{lse}{LSE}{linear system of equations}
\newacronym{msb}{MSB}{most-significant bit}
\newacronym{pde}{PDE}{partial differential equation}
\newacronym{rhs}{RHS}{right-hand side}

\DeclareMathOperator*{\argmin}{argmin}
\DeclareMathOperator{\lstCol}{maxColIdx}
\DeclareMathOperator{\fstCol}{minColIdx}
\newcommand{\tail}{\delta}

\newcommand{\compact}[1]{\acute{#1}}
\newcommand{\std}[1]{\hat{#1}}

\newcommand{\regressive}[2]{\operatorname{regressive}(+#1, #2)}
\newcommand{\progressive}[2]{\operatorname{progressive}(+#1, #2)}

\newcommand{\precU}[1]{\textcolor{pm-rb10}{#1}}
\newcommand{\precR}[1]{\textcolor{pm-rb24}{#1}}
\newcommand{\precZ}[1]{\textcolor{pm-gren}{#1}}
\newcommand{\precA}[1]{\textcolor{pm-gren}{#1}}
\newcommand{\precS}[1]{\textcolor{black}{#1}}
\newcommand{\precF}[1]{\textcolor{db-prpl}{#1}}
\newcommand{\precL}[1]{\textcolor{black}{#1}}
\newcommand{\precM}[1]{\textcolor{db-redd}{#1}}

\newcommand{\reg}{\blacktriangle}
\newcommand{\pro}{\triangledown}
\newcommand{\PRO}{\blacktriangledown}
\newcommand{\fix}{\vartriangle}
\newcommand{\ful}{\bullet}
\newcommand{\maxFulPRO}{\blacklozenge}

\newcommand{\flReg}[1]{\underline{#1}_{\reg}}
\newcommand{\flPro}[1]{\underline{#1}_{\pro}}

\newcommand{\flFix}[1]{\underline{#1}_{\fix}}
\newcommand{\flFul}[1]{\underline{#1}_{\ful}}
\newcommand{\flMaxFulPRO}[1]{\underline{#1}_{\maxFulPRO}}

\newcommand{\algosize}{\footnotesize}
\newcommand{\tablesize}{\footnotesize}

\headers{\TheTitle}{Bauer, Kohl, McCormick, and Tamstorf}

\ifpdf%
\hypersetup{%
  pdftitle={\TheTitle},
  pdfauthor={%
    Daniel Bauer, %
    Nils Kohl, %
    Stephen F. McCormick, and %
    Rasmus Tamstorf%
  }
}
\fi

\begin{document}

  \maketitle

  \begin{abstract}
    As the discretization error for the solution of a partial differential equation (PDE) decreases, the precision required to store the corresponding coefficients naturally increases.
    Storing the solution's finite element coefficients explicitly requires $\mathcal O(n \log n)$ bits of storage, where $n$ is the number of degrees of freedom (DoFs).
    This paper presents a full multigrid method to compute the solution in a compressed format that reduces the storage complexity of the solution and intermediate vectors to $\mathcal O(n)$ bits.
    This reduction allows a matrix-free implementation to solve elliptic PDEs with an overall linear space complexity.
    For problems limited by the memory capacity of current supercomputers, we expect a memory footprint reduction of about an order of magnitude compared to state-of-the-art mixed-precision methods.
    We demonstrate the applicability of our algorithm by solving two model problems.
    Depending on the PDE and polynomial degree, but irrespective of the problem size, the solution vector on the finest grid requires between \num{4} and \num{12} bits per DoF, and the residual and correction require \num{3} to \num{6} bits each.
    Additional data is stored on the coarse grids with modestly increasing bit widths toward coarser grids.
  \end{abstract}



  \section{Introduction}
  \label{sec:introduction}

  Discretization of a \gls*{pde} generally induces a \emph{discretization error} that limits the accuracy of the computed solution.
  These errors occur of course in any such setting because a solution in an infinite-dimensional space cannot be expected to be exactly represented by a finite-dimensional computation.
  In the following, our focus is on a fully regular elliptic \gls*{pde} of order $2m$ discretized by the \gls*{fem} on a uniform grid of mesh size $h$ using polynomials of degree $p$.
  In this setting, the energy norm
  of the discretization error is of order $h^{p-m+1}$~\cite{strang_analysis_2008}.

\begin{figure}
    \centering

    \subfloat{
      \begin{tikzpicture}
        \begin{groupplot}[
          group style = {group size=2 by 1, y descriptions at=edge left, horizontal sep=4mm},
          width = 0.35\textwidth,
          every axis/.append style = {
            xlabel = {FMG level $L$},
            ymode = log,
            ymax = 10,
            ymin = 10^-15,
          },
          table/x = level,
          cycle multiindex list = {[samples of colormap = 6 of viridis]\nextlist mymarklist},
        ]
          \nextgroupplot[
            title = {Single precision},
            ylabel = {Relative $H^1$-error},
            ytick = {10^0, 10^-3, 10^-6, 10^-9, 10^-12, 10^-15},
          ]
          \addplot table[y=p1sp] {data/hockey-sticks.dat};
          \addplot table[y=p2sp] {data/hockey-sticks.dat};
          \addplot table[y=p3sp] {data/hockey-sticks.dat};
          \addplot table[y=p4sp] {data/hockey-sticks.dat};
          \addplot table[y=p5sp] {data/hockey-sticks.dat};

          \begin{scope}[domain=1:15,black,dashed,mark=none]
            \addplot {0.60*0.5^(1*x)};
            \addplot {0.20*0.5^(2*x)};
            \addplot {0.10*0.5^(3*x)};
            \addplot {0.04*0.5^(4*x)};
            \addplot {0.02*0.5^(5*x)};
          \end{scope}

          \nextgroupplot[
            title = {Double precision},
            ytick = {10^0, 10^-3, 10^-6, 10^-9, 10^-12, 10^-15},
            legend to name = {legend:hockey-sticks},
            legend entries = {$p=1$, $p=2$, $p=3$, $p=4$, $p=5$, $\mathcal O(h^p)$},
          ]
          \addplot table[y=p1dp] {data/hockey-sticks.dat};
          \addplot table[y=p2dp] {data/hockey-sticks.dat};
          \addplot table[y=p3dp] {data/hockey-sticks.dat};
          \addplot table[y=p4dp] {data/hockey-sticks.dat};
          \addplot table[y=p5dp] {data/hockey-sticks.dat};

          \begin{scope}[domain=1:15,black,dashed,mark=none]
            \addplot {0.60*0.5^(1*x)};
            \addplot {0.20*0.5^(2*x)};
            \addplot {0.10*0.5^(3*x)};
            \addplot {0.04*0.5^(4*x)};
            \addplot {0.02*0.5^(5*x)};
          \end{scope}
        \end{groupplot}
      \end{tikzpicture}
    }
    \subfloat{
      \raisebox{5.5ex}{
        \ref*{legend:hockey-sticks}
      }
    }
    \caption{%
      Solving the homogeneous Poisson equation on the unit interval for the manufactured solution $u(x) = x(1-x) \cos(\frac\pi 2 x)$.
      The problem is discretized with B-spline finite elements of degree $p$ and solved by a ``textbook'' \acrlong*{fmg} solver executed in an a priori fixed floating point precision.
      The error ceases to decrease with the expected rate of $\mathcal O(h^p)$ as the limit of the finite precision is reached.
      A higher precision allows for finer grids, but any fixed bit-width constrains the attainable accuracy.%
    }
    \label{fig:hockey-sticks}
  \end{figure}
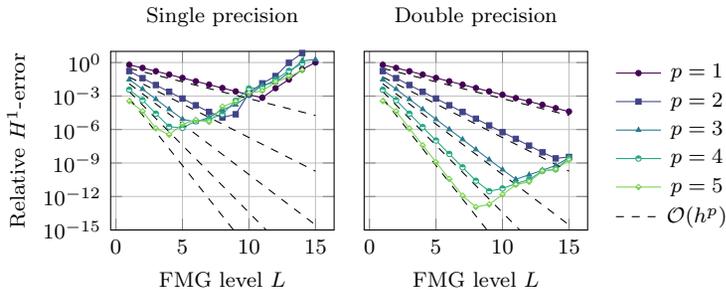

  In practice, the finite element approximation of the solution is stored as a vector of coefficients, each quantized to finite precision with unit roundoff $\varepsilon$.
  This introduces a second kind of error, the \emph{quantization error}~\cite{tamstorf_discretization-error-accurate_2021}.
  The discrete energy norm of this error is $\mathcal O(\kappa^{1/2} \varepsilon)$, where $\kappa$ is the condition number of the discretized linear operator $A$~\cite{mccormick_algebraic_2021}.
  Because $\kappa \in \mathcal O(h^{-2m})$ grows as the grid is refined~\cite{strang_analysis_2008}, so does the quantization error if $\varepsilon$ is kept fixed.
  Thus, for small enough $h$, the quantization error dominates the discretization error, leading to poor accuracy; cf.\ \cref{fig:hockey-sticks}.

  An efficient \gls*{pde} solver must aim to balance both sources of error.
  It would be a waste of computational resources to reduce one significantly, while letting the other severely limit the accuracy of the solution.
  Therefore, the precision $\varepsilon$ must be of order $h^{p+1}$, assuming that the above bounds are sharp. Given $\varepsilon$, the number of bits used to represent each coefficient of the solution is $\mathcal O(\log \varepsilon^{-1})$.
  For a regular grid, the grid spacing~$h$ is related to the number of \glspl*{dof} $n$ in that $h \approx n^{-1/d}$, where $d$ is the dimension of the domain.
  So it follows that the number of bits per coefficient must increase with $\log n$, implying that the solution vector requires $\mathcal O(n \log n)$ bits of storage.
  We emphasize that this result is about quantizing the exact solution vector and is independent of how the solution is found.

  In this work, we develop \pgls*{fmg} solver to compute the solution in a compressed format, reducing the storage complexity to $\mathcal O(n)$ bits while maintaining discretization-error-accuracy.
  In addition to the solution vector, the precisions of matrices and other vectors appearing in \gls*{fmg} must also increase as the grid is refined~\cite{mccormick_algebraic_2021,tamstorf_discretization-error-accurate_2021}.
  We show how our algorithm can be implemented with linear storage complexity overall.
  Among other things, this is achieved by operating matrix-free.

  The basic idea behind the compression is to leverage the smoothness that the solution to an elliptic \gls*{pde} exhibits by virtue of the natural assumptions about the source term and sufficient regularity.
  Due to smoothness, the solution at two neighboring locations cannot differ by more than some bounded amount, and the closer the two locations are to each other, the tighter this bound becomes.
  Loosely speaking, the difference between two \glspl*{dof} that are spatially close together must be in their trailing bits, while the information contained in the leading bits is largely redundant.

  To avoid this redundancy, the multigrid hierarchy can be used to compress the solution.
  Initially restricting the discussion to just two grids, instead of discretizing the solution on the fine grid directly, it can be split into a coarse approximation (stored on the coarse grid) and a fine correction (stored on the fine grid).
  At first this may seem counter-productive, given that it requires additionally storing the coarse-grid approximation.
  However, storing an approximation on the coarse grid is comparatively cheap, and, critically, the smoothness of the solution ensures that the magnitude of the fine-grid correction is small, meaning that it can be stored in reduced precision.
  Generalizing to a multigrid hierarchy, the solution can be disassembled into a coarsest-grid approximation and successive corrections on the finer grids.
  Due to smoothness, the corrections get progressively smaller in magnitude toward increasingly finer grids, and thus can be stored in low precision without losing accuracy.

  This idea goes back to~\cite{pan_compact_1992}, where it was termed \emph{compact multigrid data structure}.
  However, algorithms to compute the solution in compact form were presented under the assumption that the system matrix and \gls*{rhs} can be stored in fixed precision, which is not generally sufficient to achieve discretization-error accuracy~\cite{mccormick_algebraic_2021,tamstorf_discretization-error-accurate_2021}.
  Moreover, \cite{pan_compact_1992}~does not address how to compute the residual with linear space complexity, despite this computation being the most sensitive to rounding error~\cite{mccormick_algebraic_2021,tamstorf_discretization-error-accurate_2021}.
  The same compression idea as in \cite{pan_compact_1992} has been realized in \cite{thurner_ganzzahlige_1994} using hierarchical bases on sparse grids.
  However, the presented solver relies on the system matrix being diagonal, essentially eliminating its practical value.

  Compressing the output of numerical simulations has also been studied decoupled from the numerical method underlying the simulation; see~\cite{li_data_2018} for an overview.
  While this cannot reduce the memory cost of the simulation, the motivation is to reduce the size of the output on permanent storage and/or make post-processing feasible~\cite{li_data_2018}.
  The authors of~\cite{ainsworth_multilevel_2018,liang_mgard_2022} draw on the ideas of a multigrid hierarchy, and separate the input data into different scales by a series of matrix-vector products and linear system solves.
  The scales are stored on a multigrid hierarchy and then further compressed through level-dependent quantization and using an external compression algorithm.
  In~\cite{reshniak_lifting_2024,kolomenskiy_waverange_2022}, wavelets are used to achieve a similar decomposition.
  Because the input data of these methods need not be further specified, they are independent of the numerical solver and offer great flexibility.
  However, with no connection to the underlying problem/theory, it is unclear what an acceptable bound for the compression error may be.
  Our approach integrates the compression with the \gls*{fmg} solver, enabling us to select quantization precisions in such a way that even the compressed data is as accurate as the discretization allows.

  In this work, we adopt the compact representation idea of~\cite{pan_compact_1992} to develop algorithms to solve elliptic \glspl*{pde} to discretization-error accuracy, with linear bit-com\-plex\-ity in space.
  To the authors' knowledge, it is the first general implementation of a \gls*{pde} solver that can attain discretization-error accurate solutions with $\mathcal O(n)$ storage.
  For state-of-the-art problems limited by memory capacity~\cite{gmeiner_quantitative_2016}, we expect a memory footprint reduction of about an order of magnitude.

  The rest of the paper is structured as follows.
  \Cref{sec:compact-representation} describes the compact format and \cref{sec:compact-mg} introduces the compact multigrid algorithm.
  In subsequent sections, we discuss the algorithm's precision requirements (\cref{sec:precision}), optimize its memory consumption to linear bit-complexity (\cref{sec:memory-optimization}), and study the associated cost (\cref{sec:cost}).
  Finally, \cref{sec:experiments} presents numerical experiments and \cref{sec:conclusion} summarizes our conclusions.

  \section{Compact representation}
  \label{sec:compact-representation}

  To formally introduce the compact representation, we first establish some notation.
  The objective is to find the solution~$u(x)$ of an elliptic \gls*{pde}, $\mathcal{L}u(x)=f(x)$, with $u(x)$ in the solution space $V$ of functions of $x$ in the domain $\Omega$.
  To solve the \gls*{pde}, we discretize it on a set of $L+1$ nested grids corresponding to finite-dimensional subspaces, $V_0 \subset \ldots \subset V_L \subset V$.
  Each of these $n_\ell$-dimensional subspaces, $V_\ell$ for $0 \leq \ell \leq L$, is spanned by a set of finite element basis functions.
  We can collect these basis functions in a (covariant) row vector $\Phi_\ell(x)$, having $n_\ell$ elements, each element being a basis function.
  Because the spaces $V_\ell$ are nested, there are prolongation matrices $P_\ell \in \mathbb R^{n_\ell \times n_{\ell-1}}$ from level $\ell-1$ to level $\ell$ such that
  \begin{equation}
    \Phi_{\ell-1}(x) = \Phi_\ell(x) P_\ell\,.
  \end{equation}

  Using the \gls*{fem} basis, the approximate solution on level $\ell$, denoted $u_\ell(x)$, can be written as
  \begin{equation}
    u_\ell(x) = \Phi_\ell(x) \std u_\ell \in V_\ell\,,
  \end{equation}
  where $\std u_\ell \in \mathbb R^{n_\ell}$ is the (contravariant) column vector of \gls*{fem} coefficients.
  We also call $\std u_\ell$ the \emph{standard representation} of $u_\ell(x)$.

\begin{figure}
    \centering

    \pgfplotsset{
      width  = 0.38\textwidth,
      height = 0.30\textwidth,
      axis x line = middle,
      axis y line = left,
      ytick = {0},
      xtick = \empty,
      cycle multiindex list = {[samples of colormap = 4 of viridis]\nextlist mymarklist},
      enlargelimits = 0.05,
      xlabel = {$x$},
    }

    \null\hfill%
    \subfloat[Exact function.]{%
      \begin{tikzpicture}
        \begin{axis}[
          title = {$u(x)$},
        ]
          \addplot[samples=129,black] {tanh(x)+2^(-0.2*(x-1)^2)*sin(80*x)};
        \end{axis}
      \end{tikzpicture}
    }\hfill%
    \subfloat[Standard representation.]{
      \begin{tikzpicture}
        \begin{axis}[
          title = {\llap{$\approx$\hspace{5.5em}}$\Phi_2(x) \std u_2$\rlap{\hspace{4.1em}$=$}},
        ]
          \addplot+[samples=33] {tanh(x)+2^(-0.2*(x-1)^2)*sin(80*x)};
        \end{axis}
      \end{tikzpicture}
    }\hfill%
    \subfloat[Compact representation.]{
      \begin{tikzpicture}
        \begin{axis}[
          title = {$\Phi_2(x) (\compact u_2 + P_2 (\compact u_1 + P_1 \compact u_0))$},
          legend entries = {$\compact u_0$, $\compact u_1$, $\compact u_2$},
          legend style = {at={(0.03,1.03)},anchor=north west,draw=black,inner xsep=1pt},
        ]
          \addplot+[mark repeat=4] table[y=l0] {data/example-func.dat};
          \addplot+[mark repeat=2] table[y expr=\thisrow{l1}-\thisrow{l0}] {data/example-func.dat};
          \addplot+[mark repeat=1] table[y expr=\thisrow{l2}-\thisrow{l1}] {data/example-func.dat};
        \end{axis}
      \end{tikzpicture}
    }\hfill\null%

    \caption{%
      An example function discretized with piecewise-linear elements.
      The basis on level~\num{2} is denoted by the row vector $\Phi_2(x)$.
      In standard representation, only the finest grid is used.
      In contrast, the compact vector consists of a coarse approximation and two successively finer corrections.
      Crucially, the magnitude of the corrections decreases toward increasingly finer grids, allowing them to be stored in lower precision without significant loss of accuracy.%
    }
    \label{fig:compact-representation}
  \end{figure}
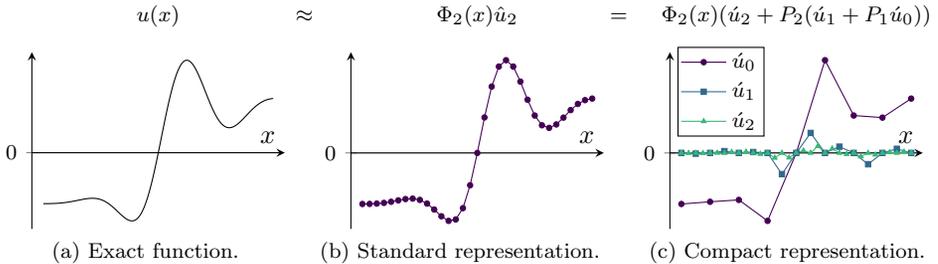

  The \emph{compact representation} of the function $u_L(x)$ is a set of $L+1$ \emph{sections} $\compact u_\ell \in \mathbb R^{n_\ell}$ for each level $\ell$, $0 \leq \ell \leq L$, which fulfills
  \begin{equation}
    \label{eq:compact-to-standard}
    \std u_L = \compact u_L + P_L (\compact u_{L-1} + P_{L-1} (\cdots + P_1 \compact u_0))\,.
  \end{equation}
  To refer to the collection of all compact sections on grids $0$ to $L$, we write $\compact u_{0\rng L}$, or $\compact u$ for brevity.
  To emphasize, we distinguish the function $u_L(x) \in V_L$, its \gls*{fem} coefficient vector $\std u_L \in \mathbb R^{n_L}$ decorated with a hat, and the coefficient vectors of the compact sections $\compact u_{0 \rng L} \in \bigotimes_{\ell=0}^L \mathbb R^{n_\ell}$ marked with an acute accent.
  A sketch of a function and its standard and compact representations is given in \cref{fig:compact-representation}.

  \subsection{Precision requirements}
  \label{sec:compact-representation-precision}

  If the hierarchy comprises at least two levels, there are infinitely many compact representations $\compact u_{0 \rng L}$ fulfilling \cref{eq:compact-to-standard} due to the nested nature of the spaces $V_\ell$~\cite{griebel_multilevel_1994}.
  Given one compact vector, any function can be added to grid $\ell$ if it is also subtracted from the next finer grid $\ell+1$ to construct a different representation of the same function.
  Even a degenerate compact representation is possible where only the finest section is nonzero.
  Thus, how many bits are required to store each section cannot be answered universally.
  In this latter example, the finest section has to equal the standard representation while all other sections are zero.
  Hence, this particular compact representation has the same storage cost as the standard representation: the precision of the finest section increases with the number of levels~\cite{mccormick_algebraic_2021,tamstorf_discretization-error-accurate_2021} while the coarse sections require zero bits.

  Storage savings are obtained by compact representations whose coarse sections closely approximate the continuous function.
  Consider the special example of a fully regular second-order elliptic \gls*{pde}.
  If the \gls*{rhs} is in $L^2(\Omega)$, then the solution is in $H^2(\Omega)$~\cite{strang_analysis_2008} and the best finite element approximations in the \emph{energy norm} (defined in \cref{lemma:smoothness} below) converge to the solution in the $H^1(\Omega)$ norm.
  The following lemma shows that the segment of the best energy approximation on a fine grid, defined as the difference between the best fine- and coarse-grid energy approximations, vanishes in the energy norm sense as $L \to \infty$.

  \begin{lemma}%
    \label{lemma:smoothness}%
    Consider a second-order \gls*{pde}, written in weak form as
    \begin{equation}
      \label{Lu=f}
      \text{Find } u \in H^1(\Omega) \text{ such that } a(u,v) = b(v) \ \forall v \in H^1(\Omega)\,,
    \end{equation}
    where $a(u,v)$ is a symmetric $H^1(\Omega)$ bounded and coercive bilinear form and $b(v)$ is a linear form, both derived from weakening the strong form $\mathcal L u = f$ of the \gls*{pde}.
    Consider the energy norm $a(\cdot,\cdot)^{1/2}$.
    As $L \to \infty$, the difference between the best energy approximations to $u$ from $V_L$ and $V_{L-1}$ vanishes in the energy norm sense.
  \end{lemma}
  \begin{proof}
    The first result here about energy orthogonality between the best energy approximation $u_L \in V_L$ to the \gls*{pde} solution $u \in H^1(\Omega)$ and its error $t := u - u_L$ is a well-known basic principle in some contexts.
    See~\cite{strang_analysis_2008}, for example.
    Note that the energy orthogonality of $u_L \in V_L$ and $t$ follows naturally from the property that $u_L$ is the energy-orthogonal projection of $u$ onto $V_L$.
    Nevertheless, the proof is included here for completeness and because an analogous argument leads to the discrete case that the same holds for the best energy approximation $u_{L-1} \in V_{L-1}$ to $u_L$.
    This discrete result then provides the basis for establishing the lemma.

    The best energy approximation $u_L \in V_L$ to $u \in H^1(\Omega)$ is, by definition,
    \begin{equation}
      \label{energynorm}
      u_L = \argmin_{v_L \in V_L} a(v_L - u, v_L - u)\,.
    \end{equation}
    To establish the energy orthogonality that we need, note first that
    \begin{equation}
      \label{ah}
      \begin{aligned}
        &a(v_L - u, v_L - u) \\
        &\qquad= a(v_L - u_L, v_L - u_L) +2 a(v_L - u_L, u_L - u) + a(u_L - u, u_L - u) \\
        &\qquad= a(v_L - u_L, v_L - u_L) + C\,,
      \end{aligned}
    \end{equation}
    where $C = a(u_L - u, u_L - u)$ is a constant w.r.t.\ $v_L$.
    The assertion here is that the term $a(v_L - u_L, u_L - u)$ vanishes because $u_L - u$ is energy-orthogonal to $V_L$.
    Were it not, then we could determine a nonzero $\delta_L \in V_L$ so that $u_L - u + \delta_L$ is energy-orthogonal to $V_L$, but that leads to the following contradiction to \cref{energynorm} for any $s \in (0, \frac{2}{3})$:
    \begin{align*}
      a(u_L &+ s \delta_L - u, u_L + s \delta_L - u) \\
      &= a(u_L - u, u_L - u) + 2s a(\delta_L, u_L + s \delta_L - u) +  s^2 a(\delta_L, \delta_L)\\
      &= a(u_L - u, u_L - u) + 2s a(\delta_L, u_L - u + \delta_L + (s - 1) \delta_L)  +  s^2 a(\delta_L, \delta_L)\\
      &= a(u_L - u, u_L - u) + 2s a(\delta_L, (s - 1) \delta_L)  +  s^2 a(\delta_L, \delta_L)\\
      &= a(u_L - u, u_L - u)  +  (3s - 2) s a(\delta_L, \delta_L) \\
      &< a(u_L - u, u_L - u)\,!
    \end{align*}

    Note now by energy orthogonality of $t$ and $V_L$ that
    \begin{align*}
      a(v_{L-1} - u, v_{L-1} - u) &= a(v_{L-1} - u_L + t, v_{L-1} - u_L + t) \\
      &= a(v_{L-1} - u_L , v_{L-1} - u_L) + a(t, t) .
    \end{align*}
    Since $a(t, t)$ is independent of $v_{L-1}$, then this expansion shows that the best energy approximation $u_{L-1}$ to $u$ from $V_{L-1}$ is the arg min of $a(v_{L-1} - u_L, v_{L-1} - u_L)$ over $v_{L-1} \in V_{L-1}$.
    An argument analogous to that above then shows that the error $t_L = u_L - u_{L-1}$ is energy-orthogonal to $V_{L-1}$.

    All of this development goes to show that $u_L = u_{L-1} + t_L$ is an en\-ergy-orthog\-o\-nal decomposition,
    which means of course that the only difference between $u_{L-1}$ and $u_L$ is the segment $t_L$.
    The energy norm of $t_L$ must therefore vanish as $L \to \infty$ because $u_L$ converges in energy to $u$, thus proving the lemma.
  \end{proof}

  By \cref{lemma:smoothness}, the segments $t_L$ of the energy-orthogonal decomposition vanish in energy as $L \to \infty$.
  This in turn implies that the  bits representing $t_L$ must be zero (relative to the full representation of $u_L$) up to a small number of trailing bits.
  This justifies the assumption in this work of the existence of a compact representation that can be stored in \emph{regressive precision}, by which we mean that
  \begin{equation}
    \label{eq:regressive}
    k (L-\ell) + b
  \end{equation}
  bits are used to store the compact section on level~$\ell$, where $k$ and $b$ are constants independent of $L$ and $\ell$.
  That is, the precision of the finest section ($b$ bits) is independent of the number of levels~$L$, and the number of bits used for the coarser sections increases (by $k$ bits) toward increasingly coarser grids.
  We will denote regressive precision as $\regressive{k}{b}$, where the plus sign serves as a reminder that the first argument is the increment.
  \Cref{sec:cost} details why the total number of bits is linear in this case.
  The term ``regressive'' is chosen in analogy to \emph{progressive} precision from~\cite{mccormick_algebraic_2021,tamstorf_discretization-error-accurate_2021}, where the precision increases with the refinement level:
  If a hierarchy of vectors is stored in $\progressive{k}{b}$ precision, $b$~bits are used on the coarsest grid, and this precision increases by $k$ bits for each level of refinement.

  To achieve discretization-error accuracy, the standard representation of the solution of the \gls*{pde} requires $\progressive{p+1}{b}$ precision~\cite{mccormick_algebraic_2021,tamstorf_discretization-error-accurate_2021}.
  Below we demonstrate experimentally that the solution's compact representation can be stored in regressive precision, using the same increment of $p+1$ bits.
  In essence, the precision hierarchy is flipped, so that now the fine-grid vectors---comprising most \glspl*{dof} of the multigrid hierarchy---are stored in fixed and low precision, instead of requiring an ever-increasing number of bits per \gls*{dof}.
  We reiterate that this requires finding a suitable compact representation with good coarse-grid approximations, which is a major concern in the development of \gls*{cmg} below.
  \cref{fig:prog-vs-reg} illustrates both ways of storing the solution.

  \begin{figure}
    \centering

    \tikzset{
      every picture/.style = {scale=0.45},
    }

    \null\hfill
    \subfloat[%
      $\progressive{6}{2}$ precision:
      Bit width increases linearly \emph{from coarse to fine} grids.
      Space complexity $\mathcal O(n_L \log n_L)$ (cf.\ \cref{sec:introduction}).%
    ]{
      \begin{tikzpicture}
        \begin{scope}[fill=pm-rb10,fill opacity=0.2,text opacity=1,font=\scriptsize]
          \begin{scope}[shift={(-3.5,3.8)}]
            \foreach \y/\z in {0/0.5, 0.8/1.3, 1.6/2.1, 2.4/2.9}{
              \filldraw ( 0.0,\y) rectangle ( 7.0,\z);
            }
            \node foreach \x/\b in {0.25/1,0.75/0,1.25/0,1.75/0,2.25/1,2.75/0,3.25/1,3.75/0,4.25/0,4.75/1,5.25/1,5.75/0,6.25/0,6.75/1} at (\x,2.65) {\b};
            \node foreach \x/\b in {0.25/1,0.75/0,1.25/0,1.75/0,2.25/1,2.75/0,3.25/1,3.75/1,4.25/0,4.75/1,5.25/0,5.75/0,6.25/0,6.75/0} at (\x,1.85) {\b};
            \node foreach \x/\b in {0.25/1,0.75/0,1.25/0,1.75/0,2.25/0,2.75/0,3.25/1,3.75/1,4.25/1,4.75/0,5.25/0,5.75/0,6.25/1,6.75/0} at (\x,1.05) {\b};
            \node foreach \x/\b in {0.25/1,0.75/0,1.25/0,1.75/0,2.25/0,2.75/0,3.25/0,3.75/0,4.25/1,4.75/0,5.25/0,5.75/1,6.25/1,6.75/0} at (\x,0.25) {\b};
          \end{scope}
        \end{scope}

        \begin{scope}[fill=pm-rb10,fill opacity=0.1,text opacity=0.8,font=\scriptsize,dash pattern=on 3pt off 1.5pt]
          \begin{scope}[shift={(-3.5,1.5)}]
            \foreach \y/\z in {0/0.5, 0.8/1.3}{
              \filldraw (0.0,\y) rectangle (4.0,\z);
            }
            \node foreach \x/\b in {0.25/1,0.75/0,1.25/0,1.75/0,2.25/1,2.75/0,3.25/0,3.75/1} at (\x,1.05) {\b};
            \node foreach \x/\b in {0.25/1,0.75/0,1.25/0,1.75/0,2.25/0,2.75/0,3.25/0,3.75/0} at (\x,0.25) {\b};
          \end{scope}

          \begin{scope}[shift={(-3.5,0)}]
            \filldraw (0.0,0.0) rectangle (1.0,0.5);
            \node foreach \x/\b in {0.25/1,0.75/0} at (\x,0.25) {\b};
          \end{scope}
        \end{scope}

        \begin{scope}[shift={(-4.5,0)}, left]
          \node[align=center] at (0,5.25) {fine grid\\(level $L$)};
          \node               at (0,0.25) {coarse grid};
        \end{scope}

        \begin{scope}[shift={(0,7.3)}]
          \draw[latex-latex] (-3.5,0) -- node[fill=white] {bit width} (3.5,0);
        \end{scope}

        \begin{scope}[dashed]
          \foreach \y in {1.0,3.3} {
            \draw (-4,\y) -- (4,\y);
          }
        \end{scope}

        \begin{scope}[fill=black,radius=3pt]
          \begin{scope}[shift={(-4.1,0)}]
            \filldraw (0,6.45) circle[] --
                      (0,5.65) circle[] --
                      (0,4.85) circle[] --
                      (0,4.05) circle[];
            \filldraw (0,2.55) circle[] --
                      (0,1.75) circle[];
            \filldraw (0,0.25) circle[];
          \end{scope}
        \end{scope}
      \end{tikzpicture}
    }\hfill
    \subfloat[%
      $\regressive{6}{2}$ precision:
      Bit width increases linearly \emph{from fine to coarse} grids.
      Space complexity $\mathcal O(n_L)$ (cf.\ \cref{sec:cost}).%
    ]{
      \begin{tikzpicture}
        \begin{scope}[fill=pm-rb10,fill opacity=0.2,text opacity=1,font=\scriptsize]
          \begin{scope}[shift={(+2.5,3.8)}]
            \foreach \y/\z in {0/0.5, 0.8/1.3, 1.6/2.1, 2.4/2.9}{
              \filldraw ( 0.0,\y) rectangle ( 1.0,\z);
            }
            \node foreach \x/\b in {0.25/0,0.75/1} at (\x,2.65) {\b};
            \node foreach \x/\b in {0.25/1,0.75/1} at (\x,1.85) {\b};
            \node foreach \x/\b in {0.25/1,0.75/1} at (\x,1.05) {\b};
            \node foreach \x/\b in {0.25/1,0.75/0} at (\x,0.25) {\b};
          \end{scope}

          \begin{scope}[shift={(-0.5,1.5)}]
            \foreach \y/\z in {0/0.5, 0.8/1.3}{
              \filldraw (0.0,\y) rectangle (4.0,\z);
            }
            \node foreach \x/\b in {0.25/0,0.75/0,1.25/1,1.75/1,2.25/1,2.75/0,3.25/1,3.75/0} at (\x,1.05) {\b};
            \node foreach \x/\b in {0.25/1,0.75/1,1.25/1,1.75/0,2.25/0,2.75/1,3.25/0,3.75/1} at (\x,0.25) {\b};
          \end{scope}

          \begin{scope}[shift={(-3.5,0)}]
            \filldraw (0.0,0.0) rectangle (7.0,0.5);
            \node foreach \x/\b in {0.25/1,0.75/0,1.25/0,1.75/0,2.25/1,2.75/0,3.25/1,3.75/0,4.25/0,4.75/1,5.25/1,5.75/0,6.25/0,6.75/1} at (\x,0.25) {\b};
          \end{scope}
        \end{scope}

        \begin{scope}[shift={(-4.5,0)}, left]
          \node[align=center] at (0,5.25) {fine grid\\(level $L$)};
          \node               at (0,0.25) {coarse grid};
        \end{scope}

        \begin{scope}[shift={(0,7.3)}]
          \draw[latex-latex] (-3.5,0) -- node[fill=white] {bit width} (3.5,0);
        \end{scope}

        \begin{scope}[dashed]
          \foreach \y in {1.0,3.3} {
            \draw (-4,\y) -- (4,\y);
          }
        \end{scope}

        \begin{scope}[fill=black,radius=3pt,shift={(-4.1,0)}]
          \filldraw (0,6.45) circle[] --
                    (0,5.65) circle[] --
                    (0,4.85) circle[] --
                    (0,4.05) circle[];
          \filldraw (0,2.55) circle[] --
                    (0,1.75) circle[];
          \filldraw (0,0.25) circle[];
        \end{scope}
      \end{tikzpicture}
    }\hfill\null

    \caption{%
      Sketch of a vector in standard representation (progressive precision) and compact representation (regressive precision).
      The level index increases from bottom to top, and bits become less significant from left to right.
      In the compact case, all sections $\compact u_{0\rng L}$ are stored as a whole, while in standard representation the finest grid vector represents the entire solution.
      The bit width of the finest compact section is independent of $L$ while the precisions of all coarser sections increase with $L$.
      The bit-values shown are made up and serve only illustrative purposes.
    }
    \label{fig:prog-vs-reg}
  \end{figure}
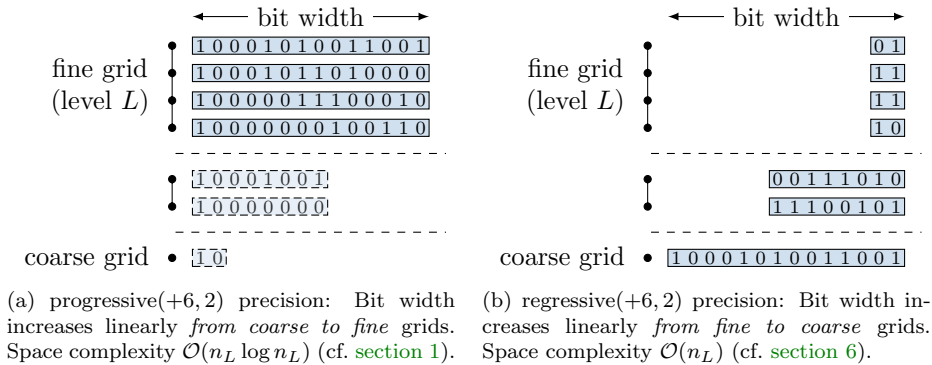

  \subsection{Compact vector operations}

  According to \cref{eq:compact-to-standard}, the \gls*{fem} coefficient vector $\std u_L$ can be recovered by summing all compact sections after prolongating them to the finest grid.
  Converting a compact representation vector to standard representation in this way is straightforward and requires only $\mathcal O(n_L)$ operations.
  The reverse conversion is more involved.
  For example, the decomposition of~\cite{ainsworth_multilevel_2018} requires obtaining $L^2$ projections onto the coarse grids, which involves solving a variational problem, and in multiresolution analysis the decomposition matrices are generally not sparse unless biorthogonal wavelets are constructed~\cite{stollnitz_wavelets_1996}.
  We will not go into details here and instead show how to compute the solution from the \gls*{pde} directly in compact form, eliminating the need for a compactification routine.

  Multiplying a compact vector with a scalar $\alpha$ is as simple as multiplying each section with $\alpha$.
  Likewise, sums and general linear combinations of compact vectors can be formed on each section individually because they commute with the prolongations in \cref{eq:compact-to-standard}.
  It is worth noting that the scalar $\alpha$ in the multiplication $\compact y_{0 \rng L} \gets \alpha_{0 \rng L} \compact x_{0 \rng L}$ can vary with the level.
  We exploited this successfully to implement an additive \gls*{bpx}-like~\cite{bramble_parallel_1990} multigrid method by varying the smoothing rate depending on the level.
  However, adequate rates are problem-specific and non-trivial to find, so we focus instead on multiplicative multigrid in this paper.

  Prolongating a compact vector amounts to introducing a zero-valued correction on the new fine-grid.
  If the vector is stored in $\regressive{k}{b}$ precision, the precision of the new finest section is $b$ bits.
  Because the number of levels $L$ increases in \cref{eq:regressive}, the precision of all other sections is incremented by appending $k$ trailing zero bits.

  \section{Compact multigrid}
  \label{sec:compact-mg}

  In this section, we discuss how the \gls*{fmg} algorithm can be adapted so that it computes the solution in compact format.
  While doing so, we focus on the storage cost of the compact solution and address computing a compact representation that can be stored in regressive precision, i.e., with linear space complexity.
  In contrast, we defer reducing the memory footprint of other quantities to \cref{sec:memory-optimization}.
  Second, we aim at keeping the excellent convergence properties of \gls*{fmg}, i.e., linear complexity w.r.t.\ the number of arithmetic operations.
  Note that here and throughout the paper we count all \glspl*{flop} equally, although their precisions influence their costs in terms of execution time, energy consumption, and chip area.
  While a standard \gls*{fmg} implementation performs all \glspl*{flop} in the same bit-width, which increases with $n$, our algorithm reduces the precision of some of these \glspl*{flop} and never requires more precision than the standard implementation.

  We base our algorithm on the \gls*{fas} variant of multigrid, briefly defined below.
  For details, the reader is referred to~\cite{brandt_multigrid_2011,briggs_multigrid_2000}.
  For convenience, the supplement distills the relevant ideas into this article's notation.
  \Gls*{fas} owes its name to the fact that it computes a full approximation to the fine-grid solution on each grid of the hierarchy.
  Precisely this property is helpful when computing a compact solution, so we choose \gls*{fas} over the equivalent and more ``standard'' \gls*{cs}.

  \subsection{Full approximation scheme}
  \label{sec:fas}

  The objective is to solve an elliptic \gls*{pde} discretized on refinement level $L$ as
  \begin{equation}
    \label{eq:fine-grid-std}
    A_L \std y_L = \std r_L\,.
  \end{equation}
  We further define standard multigrid components~\cite{brandt_multigrid_2011,briggs_multigrid_2000}:
  \begin{itemize}
    \item the prolongation matrices $P_\ell \in \mathbb R^{n_\ell \times n_{\ell-1}}$, $1 \leq \ell \leq L$, that canonically arise from the nesting of the \gls*{fem} spaces as outlined in the previous section,
    \item the restriction matrices $R_\ell \in \mathbb R^{n_{\ell-1} \times n_\ell}$, $1 \leq \ell \leq L$, that transfer vectors from one grid to the next coarser grid and are typically defined as $R_\ell := P_\ell^T$,
    \item coarse-grid operators $A_\ell := R_{\ell+1} A_{\ell+1} P_{\ell+1}$, $0 \leq \ell < L$, arising through Ga\-ler\-kin coarsening,
    \item and, for all $\ell$, $0 \leq \ell \leq L$, a smoother of the form~\cite{falgout_generalizing_2004,baker_multigrid_2011}
      \begin{equation}
        \label{eq:smoother}
        \std y_\ell \gets \std y_\ell + M_\ell (\std r_\ell - A_\ell \std y_\ell)\,.
      \end{equation}
  \end{itemize}

  Both \gls*{cs} and \gls*{fas} use the smoother to reduce oscillatory errors and obtain an approximate solution $\std y_L$.
  While \gls*{cs} uses the next coarser grid to find a smooth correction $\std d_{L-1}$, the \gls*{fas} coarse-grid variable approximates the full solution to the fine-grid equation.
  It is defined as
  \begin{equation}
    \label{eq:coarse-grid-variable-fas}
    \std y_{L-1} := \tilde R_L \std y_L + \std d_{L-1}\,,
  \end{equation}
  where $\tilde R_L$ is a restriction operator defined on the finite-element solution space $V_L$, unlike $R_L$ which is defined on the dual space of $V_L$.

  In \gls*{fas}, $\std y_{L-1}$ is found by (approximately) solving the coarse-grid equation
  \begin{equation}
    \label{eq:coarse-grid-fas}
      A_{L-1} \std y_{L-1} = R_L \std r_L + \underbrace{(A_{L-1} \tilde R_L \std y_L - R_L A_L \std y_L)}_{\std \tau_L}\,,
  \end{equation}
  where $\std \tau_L$ is known as the \emph{fine-to-coarse defect correction}.
  It modifies the coarse-grid equation so that at convergence its solution coincides with the (restricted) fine-grid solution $\tilde R_L A_L^{-1} \std r_L$.

  After having found a good approximation on the coarse grid through recursion, the fine-grid variable is updated according to
  \begin{equation}
    \label{eq:fas-apply-correction}
    \std y_L \gets \std y_L + P_L \std d_{L-1} = \std y_L + P_L (\std y_{L-1} - \tilde R_L \std y_L)\,.
  \end{equation}

  The key difference between \gls*{cs} and \gls*{fas} is the addition of $A_{L-1} \tilde R_L \std y_L$ on both sides of \cref{eq:coarse-grid-fas}.
  Therefore, in absence of rounding errors, \gls*{cs} and \gls*{fas} have the same convergence behavior and compute exactly the same \emph{fine-grid} solution for any choice of $\tilde R_L$~\cite{brandt_multigrid_2011}.
  Moreover, with the particular choice of $\tilde R_L := 0$, \gls*{fas} reduces to \gls*{cs}, and both schemes compute the same \emph{coarse-grid} solutions.
  In essence, \gls*{fas} introduces the primal-space restriction $\tilde R_L$ that can be chosen freely to define the coarse-grid variables.
  Note that $\tilde R_L$ is independent of the dual-space restriction $R_L$.
  Remember from \cref{sec:compact-representation-precision} that the compact representation leads to storage savings only if the coarse sections represent the continuous function well.
  Hence, in the compact context, $\tilde R_L$ is a key component that must be chosen so that the coarse-grid variable $\std y_{L-1}$ closely approximates the fine-grid variable $\std y_L$.

  \subsection{Compact full approximation scheme}
  \label{sec:compact-fas}

  Now we adapt \gls*{fas} to compute a compact solution $\std y_L = \compact y_L + P_L (\compact y_{L-1} + P_{L-1} (\cdots + P_1 \compact y_0))$.
  Inserting the compact representation into the fine-grid equation~\cref{eq:fine-grid-std} gives
  \begin{equation}
    \label{eq:residual-compact-y}
    A_L (\compact y_L + P_L (\compact y_{L-1} + P_{L-1} (\cdots + P_1 \compact y_0))) = \std r_L\,.
  \end{equation}
  Before applying a smoother to \cref{eq:residual-compact-y}, the equation must be reformulated w.r.t.\ a single solution variable.
  To that end, note that the smoother primarily affects oscillatory components and, hence, its main effect should be on the finest compact section.
  This motivates moving the coarse sections to the \gls*{rhs}, giving an equation for just the finest section:
  \begin{equation}
    \label{eq:fine-grid-compact}
    A_L \compact y_L = \std r_L - \underbrace{A_L (P_L (\compact y_{L-1} + P_{L-1} (\cdots + P_1 \compact y_0)))}_{\std v_L}\,.
  \end{equation}
  If a smoother of the form \cref{eq:smoother} is used, then smoothing on \cref{eq:fine-grid-compact} is equivalent to smoothing on the original fine-grid equation~\cref{eq:fine-grid-std}, i.e., converting $\compact y_{0\rng L}$ to standard representation via \cref{eq:compact-to-standard} after applying the smoother to \cref{eq:fine-grid-compact} gives the same result as smoothing on \cref{eq:fine-grid-std}.
  To see this, insert \cref{eq:fine-grid-compact} into \cref{eq:smoother}:
  \begin{equation}
    \compact y_L \gets \compact y_L + M_L ( \std r_L - \std v_L - A_L \compact y_L )
                   =   \compact y_L + M_L ( \std r_L - A_L \std y_L )\,.
  \end{equation}
  Comparing to \cref{eq:smoother} establishes that the solution variable is changed by $M_L ( \std r_L - A_L \std y_L )$ irrespective of whether the smoother is applied to \cref{eq:fine-grid-std} or \cref{eq:fine-grid-compact}.

  After the smoothing step, we turn our attention to solving the \gls*{fas} coarse-grid equation~\cref{eq:coarse-grid-fas}.
  This requires defining the primal-space restriction operator $\tilde R$.
  Because this operator is only applied to the compact vector $\compact y_{0\rng L}$, the canonical choice is to drop the finest section:
  \begin{equation}
    \tilde R_L (\compact y_L + P_L (\compact y_{L-1} + P_{L-1} (\cdots + P_1 \compact y_0))) := \compact y_{L-1} + P_{L-1} (\cdots + P_1 \compact y_0)\,.
  \end{equation}
  Implicitly defining the restriction in this way results in good coarse approximations that allow for regressive precision storage, i.e., storage savings.
  This is because $\tilde R_L$ only removes the finest section $\compact y_L$, which stores primarily oscillating components.
  Furthermore, due to this definition (and the Galerkin condition $A_{L-1} = R_L A_L P_L$), the defect correction from \cref{eq:coarse-grid-fas} simplifies to $\std \tau_L = -R_L A_L \compact y_L$.
  Remarkably, $\std \tau_L$ depends only on the finest section $\compact y_L$ instead of on the entire vector $\std y_L$ as in~\cref{eq:coarse-grid-fas}.
  Again, moving all coarser sections to the \gls*{rhs} gives an equation for the section $\compact y_{L-1}$:%
  \begin{equation}
      \label{eq:coarse-grid-compact}
      A_{L-1} \compact y_{L-1} = R_L \std r_L - \underbrace{R_L A_L \compact y_L}_{-\std \tau_L} - \underbrace{A_{L-1} (P_{L-1} \compact y_{L-2} + \cdots + P_{L-1} \cdots P_1 \compact y_0)}_{\std v_{L-1}}\,.
  \end{equation}
  Like the fine-grid equation \cref{eq:fine-grid-compact}, the coarse-grid equation \cref{eq:coarse-grid-compact} is solved by applying a smoother and solving a coarse-grid equation.
  To display the pattern, we give the equation on the next coarser grid, which involves another $\std \tau$-term but otherwise looks familiar:
  \begin{equation}
    \label{eq:next-coarser-grid-compact}
    A_{L-2} \compact y_{L-2} = R_{L-1} R_L \std r_L - R_{L-1} \underbrace{R_L A_L \compact y_L}_{-\std \tau_L} - \underbrace{R_{L-1} A_{L-1} \compact y_{L-1}}_{-\std \tau_{L-1}} - \std v_{L-2}\,.
  \end{equation}
  This process is repeated until a coarse enough grid is reached, where solving the equation is trivial.

  To summarize, compact \gls*{fas} improves the approximation $\compact y_{0\rng L}$ by relaxing for all levels $\ell = 0,\dots,L$ on the equation
  \begin{equation}
    \label{eq:cfas}
    \begin{aligned}
      A_\ell \compact y_\ell ={}& \std r_\ell - \std s_\ell - \std v_\ell\,,\\
      \std r_\ell :={}& R_{\ell+1} \cdots R_L \std r_L\,,\\
      \std s_\ell :={}& -(\std \tau_{\ell+1} + \cdots + \std \tau_L) = R_{\ell+1} (A_{\ell+1} \compact y_{\ell+1} + R_{\ell+2} (\cdots + A_L \compact y_L))\,,\\
      \std v_\ell :={}& A_\ell \std z_\ell\,,\\
      \std z_\ell :={}& P_\ell (\compact y_{\ell-1} + P_{\ell-1} (\cdots + P_1 \compact y_0))\,.
    \end{aligned}
  \end{equation}
  The order in which the levels are handled can be chosen freely to determine the cycle structure.
  For example, a V(1,0)-cycle handles the finest level first and the coarsest level last, while a V(0,1)-cycle works from the coarsest toward the finest grid.
  Unlike in a standard \gls*{cs} or \gls*{fas} cycle, no coarse-grid corrections are applied, because smooth and oscillatory components are separated in the compact representation.

  \begin{algorithm}[tbp]
    \algosize
    \caption{\pr{cFas}: Compact V(0,1) \acrlong*{fas}.}
    \label{alg:cfas-naive}%

    \paragraph{Input}
    \begin{tabular}{ll}
      \precA{$A_{0\rng L}$}          & matrix hierarchy \\
      \precA{$M_{0\rng L}$}          & smoother \\
      \precA{$P_{1\rng L}$}          & prolongation matrices \\
      \precA{$R_{1\rng L}$}          & restriction matrices \\
      \precR{$\compact y_{0\rng L}$} & initial guess of solution \\
      \precR{$\std r_{0\rng L}$}     & \gls*{rhs} \\
    \end{tabular}

    \paragraph{Output}
    Solution $\precR{\compact y_{0\rng L}}$ improved by one V(0,1) multigrid cycle.
    \unskip\medskip

    \begin{pseudo}[kw,fullwidth,setup-append={\pseudoanchor{L-\arabic*}}]*&
      \ct{Fine-to-coarse sweep.} \\
      $\precS{\std s_L \gets 0}$ \label{line:cfas-n-f2c-start} \\
      for $\ell = (L-1) \dts 0$ \kw{step} $-1$ \\+
        $\precS{\std s_\ell \gets \precA{R_{\ell+1} A_{\ell+1} \flPro{\precR{\compact y_{\ell+1}}} + R_{\ell+1} \precS{\std s_{\ell+1}}}}$ \label{line:cfas-n-f2c-end} \\-*&

      \ct{Coarse-to-fine sweep.} \\
      $\precZ{\std z_0 \gets 0}$ \label{line:cfas-n-z0} \\
      $\precR{\compact y_0 \gets \flFix{\precM{\flMaxFulPRO{\precR{\compact y_0}} + \flMaxFulPRO{\precA{M_0}} (\flMaxFulPRO{\precR{\std r_0}} - \flMaxFulPRO{\precS{\std s_0}} - \flMaxFulPRO{\precA{A_0}} \flMaxFulPRO{\precR{\compact y_0}})}}}$ \hfill\ct{Smooth $A_0 \compact y_0 = \std r_0 - \std s_0$} \label{line:cfas-n-y0} \\
      for $\ell = 1\dts L$ \kw{step} $1$ \label{line:cfas-c2f-loop} \\+
        $\precZ{\std z_\ell \gets \flPro{\precM{\flMaxFulPRO{\precA{P_\ell}} \flMaxFulPRO{\precR{\compact y_{\ell-1}}} + \flMaxFulPRO{\precA{P_\ell}} \flMaxFulPRO{\precZ{\std z_{\ell-1}}}}}}$ \label{line:cfas-n-zl} \\
        $\precR{\compact y_\ell \gets \flFix{\precM{\flMaxFulPRO{\precR{\compact y_\ell}} + \flMaxFulPRO{\precA{M_\ell}} (\flMaxFulPRO{\precR{\std r_\ell}} - \flMaxFulPRO{\precS{\std s_\ell}} - \flMaxFulPRO{\precA{A_\ell}} \flMaxFulPRO{\precZ{\std z_\ell}} - \flMaxFulPRO{\precA{A_\ell}} \flMaxFulPRO{\precR{\compact y_\ell}})}}}$ \hfill\ct{Smooth $A_\ell \compact y_\ell = \std r_\ell - \std s_\ell - A_\ell \std z_\ell$} \label{line:cfas-n-yl} \\-

      return \precR{$\compact y_{0\rng L}$}
    \end{pseudo}
    \tikz[overlay, remember picture, decoration=brace] {
      \draw[decorate] ($(L-1.north)+(18em,2.5ex)$) -- node[anchor=west,xshift=0.5em,text width=22em] {\ctfont{%
        If the initial guess is zero ($\compact y = 0$), then $\std s$ is zero, requires no storage, and the fine-to-coarse sweep can be left out.
        In \cref{sec:cfmg}, we wrap this algorithm in \acrlong*{ir} so that the initial guess is indeed zero.
      }} ($(L-3.south)+(18em,2.5ex)$);
    }
  \end{algorithm}

  \Cref{alg:cfas-naive} lists the steps of the compact V(0,1) \gls*{fas} scheme though the algorithm can easily be modified to include more pre-/post-smoothing steps.
  Later, in \cref{sec:precision}, we discuss how much precision is needed for the various quantities appearing in the algorithm.
  While these results should not concern us yet, different precisions are represented by different colors and change of precision is denoted by underlines.

  \Cref{alg:cfas-naive} starts with a sweep from finest to coarsest grid, computing $\std s$ on all grids (\crefrange{line:cfas-n-f2c-start}{line:cfas-n-f2c-end}).
  Importantly, the recursive structure of $\std s$ is exploited to keep the number of \glspl*{flop} linear.
  The fine-to-coarse sweep is followed by a coarse-to-fine sweep that computes the prolongated standard representation $\std z_\ell := P_\ell (\compact y_{\ell-1} + P_{\ell-1} (\cdots + P_1 \compact y_0))$ (\cref{line:cfas-n-z0,line:cfas-n-zl}) and applies the smoother to \cref{eq:cfas} (\cref{line:cfas-n-y0,line:cfas-n-yl}) on every grid.
  Again, we benefit from the recursive structure of $\std z$, which allows us to maintain a linear number of \glspl*{flop}.

  Another way to present the compact \gls*{fas} equations \cref{eq:cfas} is as a block-system
  \begin{equation}%
    \label{eq:block-system}
    \begin{gathered}
      \mathcal A_L
      \begin{pmatrix}
        \compact y_L \\ \compact y_{L-1} \\ \vdots \\ \compact y_0
      \end{pmatrix}
      =
      \begin{pmatrix}
        \std r_L \\ R_L \std r_L \\ \vdots \\ R_1 \cdots R_L \std r_L
      \end{pmatrix}\,,\\
      \mathcal A_L :=
      \begin{pmatrix}
        A_L & A_L P_L & \cdots & A_L P_L \cdots P_1 \\
        R_L A_L & A_{L-1} & \cdots & A_{L-1} P_{L-1} \cdots P_1 \\
        \vdots &  & \ddots & \vdots \\
        R_1 \cdots R_L A_L & R_1 \cdots R_{L-1} A_{L-1} & \cdots & A_0
      \end{pmatrix}\,.
    \end{gathered}
  \end{equation}
  Here, the terms $\std s$ and $\std v$ have been moved to the left-hand side leading to all the off-diagonal blocks.
  Multiplication of the block-lower-triangular part of $\mathcal A_L$ and $\compact y$ gives $\std s$, and by multiplying the block-upper-triangular part with $\compact y$, we get $\std v$.

  Remarkably, a Gauss-Seidel iteration on~\cref{eq:block-system} is equivalent to a V-cycle with one Gauss-Seidel pre-smoothing step~\cite{griebel_multilevel_1994}.
  Similarly, when using Gauss-Seidel in \Cref{alg:cfas-naive}, it is equivalent to Gauss-Seidel on~\cref{eq:block-system}, when the blocks are handled in reverse order.
  \Cref{eq:block-system} has previously been derived in~\cite{griebel_multilevel_1994} by discretizing the linear operator from the full set of basis functions from all levels.
  There, it is noted that this system is semi-definite, and various multilevel methods can be implemented by an iterative method on the block-system.
  For example, an additive V-cycle can be implemented by applying block-Jacobi~\cite{griebel_multilevel_1994}.

  \subsection{Residual computation}

  Although \cref{alg:cfas-naive} can be used as a standalone solver, it must be embedded in \pgls*{fmg} scheme to be sure to achieve dis\-cret\-iza\-tion-error accuracy with linear work.
  Moreover, we want to benefit from the analysis of~\cite{mccormick_algebraic_2021,tamstorf_discretization-error-accurate_2021} and wrap the V-cycle into an \gls*{ir} loop.
  Thus, we need a way to determine the residual $\std r = \std f - A \std u$, where $\std u$ is an approximate solution to the equation $A \std u = \std f$.
  The difficulty is that only the compact representation $\compact u$ is available.

  Our proposed algorithm (\Cref{alg:fmg-residual-naive}) to compute the residual from the compact solution $\compact u_{0 \rng L}$ is straightforward:
  First, the standard representation $\std u_L$ corresponding to $\compact u_{0 \rng L}$ is computed in a sweep from coarsest to finest grid by evaluating \cref{eq:compact-to-standard} (\crefrange{line:decode-start}{line:decode-end}).
  Next, the actual residual computation follows on \cref{line:compute-r}.
  After the computation, the result can be truncated to a potentially lower precision (\cref{line:truncate-r}).
  Finally, a fine-to-coarse sweep through the grid hierarchy determines the restriction of the residual on all coarser grids (\crefrange{line:restrict-start}{line:restrict-end}).
  This is done so that the hierarchy of residuals $\std r_{0\rng L}$ can serve as input to the V-cycle \cref{alg:cfas-naive}.
  We remark that the structure of \cref{alg:fmg-residual-naive} resembles an up-side-down V-cycle.

  Because the precision required to store the standard representation $\std u_\ell$ grows with the level $\ell$~\cite{mccormick_algebraic_2021}, implementing \cref{alg:fmg-residual-naive} with linear space complexity requires some care.
  This issue is discussed in \cref{sec:memory-optimization}.

  \begin{algorithm}[tbp]
    \algosize
    \caption{\pr{fmgResidual}: \gls*{fmg} residual computation.}
    \label{alg:fmg-residual-naive}%

    \paragraph{Input}
    \begin{tabular}{ll}
      \precF{$A_L$}                   & matrix \\
      \precF{$P_{1\rng L}$}           & prolongation matrices \\
      \precF{$R_{1\rng L}$}           & restriction matrices \\
      \precF{$\std f_L$}              & \gls*{rhs} \\
      \precU{$\compact u_{0 \rng L}$} & solution approximation \\
    \end{tabular}

    \paragraph{Output}
    Residual $\precR{\std r_\ell} = R_{\ell+1} \cdots R_L (\std f_L - A_L (\compact u_L + \cdots + P_L \cdots P_1 \compact u_0))$ for $0 \leq \ell \leq L$.
    \unskip\medskip

    \begin{pseudo}[kw]*&
      \ct{Decode solution.} \\
      $\std u_0 \gets \precU{\compact u_0}$ \label{line:decode-start} \\
      for $\ell = 1 \dts L$ \kw{step} $1$ \\+
        $\std u_\ell \gets \flFul{\precU{\compact u_\ell}} + \flFul{\precM{\flMaxFulPRO{\precF{P_\ell}}} \precM{\flMaxFulPRO{\precL{\std u_{\ell-1}}}}}$ \label{line:decode-end} \\-*&

      \ct{Compute residual.} \\
      $\std t_L \gets \flFul{\precM{\flMaxFulPRO{\precF{\std f_L}} - \flMaxFulPRO{\precF{A_L}} \flMaxFulPRO{\precL{\std u_L}}}}$ \label{line:compute-r} \\
      $\precR{\std r_L \gets \flFix{\precL{\std t_L}}}$ \label{line:truncate-r} \\*&

      \ct{Restrict residual.} \\
      for $\ell = (L-1) \dts 0$ \kw{step} $-1$ \label{line:restrict-start} \\+
        $\std t_\ell \gets \flFul{\precM{\flMaxFulPRO{\precF{R_{\ell+1}}}} \precM{\flMaxFulPRO{\precL{\std t_{\ell+1}}}}}$ \\
        $\precR{\std r_\ell \gets \flFix{\precL{\std t_\ell}}}$ \label{line:restrict-end} \\-

      return \precR{$\std r_{0 \rng L}$}
    \end{pseudo}
  \end{algorithm}

  \subsection{Compact full multigrid}
  \label{sec:cfmg}

  Finally, we have all the pieces to present the compact \gls*{fmg} solver for the equation $\mathcal L u(x) = f(x)$ in \cref{alg:cfmg}.
  In \cref{line:disc0,line:discL}, the linear operator $\mathcal L$ and \gls*{rhs} $f(x)$ are discretized on level $L$ before forming the matrix hierarchy by Galerkin coarsening.
  The coarse-grid solver in \cref{line:coarse-grid-solve} is exactly the same as in a standard method---on the coarsest grid, the standard and compact representation coincide.

  Every time a new \gls*{fmg} level is reached, the current solution approximation is prolongated to this level in \cref{line:prolongation}.
  Because it is represented compactly, this amounts to adding a new zero-initialized grid; no matrix-vector products need to be evaluated.
  Remember from \cref{eq:regressive} that the number of bits of a regressive-precision quantity on level~$\ell$ depends on the total number of levels~$L$.
  Therefore, incrementing $L$ and adding a new grid also involves increasing the precision of all other levels, so that regressive precision is maintained throughout the hierarchy.

  Next, the residual is computed (\cref{line:cfmg:compute-r}) and subsequently used to solve---using a few V-cycles---for the correction (\crefrange{line:cfmg:init-y}{line:cfmg:compute-y}).
  Finally, the correction is applied to the solution in \cref{line:correction} by adding each section of the compact correction to the solution.
  Residual computation and correction are repeated for a sufficient number of \gls*{ir} iterations $N$.

  \begin{algorithm}[tbp]
    \algosize
    \caption{\pr{cFmg}: Compact \gls*{fmg}.}
    \label{alg:cfmg}

    \paragraph{Input}
    \begin{tabular}{ll}
      $\mathcal L$       & continuous linear operator \\
      $f(x)$             & continuous \acrshort*{rhs} \\
      $L_\text{max}$     & index of finest grid \\
      $N$                & number of \acrshort*{ir} iterations per \gls*{fmg} level \\
      $\nu$              & number of V-cycles per \acrshort*{ir} iteration \\
    \end{tabular}

    \paragraph{Output}
    Compact solution $\compact u_{0\rng L_\text{max}}$.
    \unskip\medskip

    \begin{pseudo}[fullwidth]
      $(\precF{A_0}, \precF{\std f_0}) \gets \pr{discretize}(\mathcal L, f(x), L=0)$ \label{line:disc0} \\
      $\precU{\compact u_0 \gets \flReg{\precF{A_0^{-1} \std f_0}}}$ \label{line:coarse-grid-solve} \hfill\ct{Coarse grid solve.} \\
      \kw{for} $L=1 \dts L_\text{max}$ \kw{step} $1$ \\+
        \precU{$\compact u_L \gets 0$} \hfill\ct{Prolongate $\compact u_{0\rng L-1}$.} \label{line:prolongation} \\
        $(\precF{A_{0\rng L}}, \precF{M_{0\rng L}}, \precF{P_{1\rng L}}, \precF{R_{1\rng L}}, \precF{\std f_L}) \gets \pr{discretize}(\mathcal L, f(x), L)$ \label{line:discL} \\

        \kw{repeat} $N$ \kw{times} \\+
          $\precR{\std r_{0\rng L} \gets{}} \pr{fmgResidual}(\precF{A_L}, \precF{P_{1\rng L}}, \precF{R_{1\rng L}}, \precF{\std f_L}, \precU{\compact u_{0\rng L}})$ \label{line:cfmg:compute-r} \hfill\ct{\cref{alg:fmg-residual-naive}.} \\

          \precR{$\compact y_{0\rng L} \gets 0$} \label{line:cfmg:init-y} \\
          \kw{repeat} $\nu$ \kw{times} \label{line:nu} \\+
            $\precR{\compact y_{0\rng L} \gets{}} \pr{cFas}(\precA{\flPro{\precF{A_{0\rng L}}}}, \precA{\flPro{\precF{M_{0\rng L}}}}, \precA{\flPro{\precF{P_{1\rng L}}}}, \precA{\flPro{\precF{R_{1\rng L}}}}, \precR{\compact y_{0\rng L}}, \precR{\std r_{0\rng L}})$ \label{line:cfmg:compute-y} \hfill\ct{V-cycle, \cref{alg:cfas-naive}.} \\-

          $\precU{\compact u_{0\rng L} \gets \compact u_{0\rng L} + \flReg{\precR{\compact y_{0\rng L}}}}$ \label{line:correction} \hfill\ct{Apply correction.} \\---
      \kw{return} $\compact u_{0\rng L_\text{max}}$
    \end{pseudo}
  \end{algorithm}

  \section{Precision and rounding}
  \label{sec:precision}

  While the algorithms presented in the previous section can be used ``as is'', the sole motivation for formulating them is to put strict limits on the precision of vectors, with the goal of reducing the total memory footprint to $\mathcal O(n)$ bits.
  Therefore, we now focus our attention on which precision is required for vector and matrix entries in the above algorithms.
  Some of these precisions are based on the theory of \cite{mccormick_algebraic_2021,tamstorf_discretization-error-accurate_2021}, while everything not covered by existing theory is based on experimental evidence.
  To ensure that as many results from \cite{mccormick_algebraic_2021,tamstorf_discretization-error-accurate_2021} can be carried over, in the remaining parts of the paper, we restrict the number of V-cycles per \gls*{ir} iteration to one, i.e., we set $\nu = 1$ in \cref{alg:cfmg}.
  A convenient consequence is that the initial guess $\compact y_{0\rng L}$ in the V-cycle \cref{alg:cfas-naive} is zero.
  Thus, the initial fine-to-coarse sweep need not be computed, and $\std s_\ell = 0$ for all $\ell$ and need not be stored.
  Hence, we do not specify a precision for $\std s$.
  However, we remark that in practice it might be more efficient to increase the number of V-cycles per \gls*{ir} iteration to reduce the number of residual calculations, even though it may come at the cost of requiring slightly higher precision in the inner solver.

  The precisions that we use in our algorithms are listed in \cref{tab:precisions}.
  The values of $b_1, \ldots, b_4$ are problem-specific and are discussed in \cref{sec:experiments}.
  We require that the matrices $A$, $P$, and $R$ and \gls*{rhs} $\std f$ are free of any error besides that stemming from quantization.
  For example, if quadrature is used to evaluate integrals, the resulting error must be insignificant compared to quantization.
  Regressive and progressive precisions change over the course of the \gls*{fmg}-cycle (w.r.t.\ $L$) and during each V-cycle (w.r.t.\ $\ell$), cf.~\cref{eq:regressive}.%
  \begin{table}[tbp]
    \tablesize
    \centering
    \caption{Precisions used for vectors and matrix entries.
      $L$ is the current \gls*{fmg} level.
      The constants $b_1, \dots, b_4$ depend on the smoothness of the solution and \gls*{rhs}.
      We experimentally determine suitable constants for some example problems in \cref{sec:experiments}.
    }
    \label{tab:precisions}

    \begin{tabular}{ccc}
      \toprule
      Quantity & Precision & Symbol \\
      \midrule
      $\compact u$ & \precU{$\regressive{p+1}{b_1}$} & $\precU\reg$ \\
      $\std r$, $\compact y$ & \precR{$\regressive{m}{b_2}$} & $\precR\fix$ \\
      matrices, $\std f$ in \pr{fmgResidual} & \precF{$\progressive{p+m+1}{b_3}$} & $\precF\PRO$ \\
      matrices, $\std z$ in \pr{cFas} & \precA{$\progressive{m}{b_4}$} & $\precA\pro$ \\
      $\std u$, $\std t$ in \pr{fmgResidual} & \precL{$(p+1)L+b_1$ bits} & $\precL\ful$ \\
      \bottomrule
    \end{tabular}
  \end{table}

  To aid the reader, quantities in \crefrange{alg:cfas-naive}{alg:cfmg} are colored according to their precision listed in \cref{tab:precisions}.
  Frequently, quantities stored in different precisions are used together in an operation, e.g., a vector addition.
  To make it clear which precision the arithmetic is performed in and what the precision of the result is, we explicitly adapt the precision of (one of) the input argument(s) so that both operands share the same precision.
  We implement decreasing the precision as rounding toward the nearest representation in the lower precision, breaking ties to even, and increasing the precision as extending the binary representation by additional zero bits.
  To keep the linear algebra in the pseudocodes comprehensible, we denote change of precision by an underline and a symbol.
  The symbol for each precision is listed in \cref{tab:precisions}.
  Once the operands have been converted, the operation is performed within this precision, subject to the usual \gls*{fp} rounding behavior.
  In case of matrix vector products / stencil applications, this entails several rounding steps internally.
  For example, $A_\ell \flFix{\compact u_\ell}$ indicates that the precision of the vector is adapted before evaluating the product in the precision of the matrix.
  By contrast, in $\flFix{A_\ell \compact u_\ell}$, matrix and vector already share the same precision and the result is converted after computing the matrix-vector product.

  As discussed previously, the solution $\compact u$ is stored in regressive precision.
  It uses the same precision increment as dictated by existing theory~\cite{mccormick_algebraic_2021,tamstorf_discretization-error-accurate_2021}, but crucially the increment is now toward coarser grids.
  The constant $b_1$, i.e., the number of bits used on level $L$, depends on the \gls*{pde}, the discretization, etc.
  Likewise, the correction and residual can be stored in regressive precision.
  As before, we adopt the same increment as used in~\cite{mccormick_algebraic_2021,tamstorf_discretization-error-accurate_2021}, which is $m$ bits,
  although it is possible that better choices exist for some \glspl*{pde}.
  Since the representation of matrices and \gls*{rhs} remains unchanged compared to previous methods, we cannot expect any storage improvements.
  Hence, we use the same precision as in~\cite{mccormick_algebraic_2021,tamstorf_discretization-error-accurate_2021}.
  Finally, the temporary vectors $\std z$, $\std u$, and $\std t$ correspond to the correction in the V-cycle of~\cite{tamstorf_discretization-error-accurate_2021}, the standard representation of the solution, and the residual in the \gls*{fmg}-cycle of~\cite{tamstorf_discretization-error-accurate_2021}, respectively.
  Therefore, they use the same respective precisions as in~\cite{tamstorf_discretization-error-accurate_2021}.
  For $\std u$ and $\std t$, this is the number of bits used for the solution's coarsest compact section.
  Note that this precision is progressive in $L$ but does not depend on the level index $\ell$ in \cref{alg:fmg-residual-naive}.

  In both algorithms, progressive precision matrices are multiplied with vectors that do not use progressive precision.
  On the coarse grids, the number of bits in the vector is much higher than in the matrix, whereas on the fine grids, the matrix uses more bits.
  To ensure sufficient accuracy of matrix-vector products, we extend the precision of the arguments before performing the operation (arithmetic and intermediate results).
  Concretely, we take the precision of the respective matrix on the given level or the precision of the vector's coarsest section ($\compact y_0$ in \cref{alg:cfas-naive}, $\compact u_0$ in \cref{alg:fmg-residual-naive}), whichever is higher.
  In either case, the precision is denoted by red font and conversion to that precision by $\precM{\maxFulPRO}$.
  The number of bits depends on $\ell$ and $L$, and is different in both algorithms.
  We do not claim that this choice is necessarily optimal.

  \section{Memory optimization}
  \label{sec:memory-optimization}

  As shown in \cref{tab:precisions}, the solution, residual, and correction can be stored in regressive precision with linear storage cost.
  However, the matrices $A$, $P$, and $R$, the \gls*{rhs}, and some temporary vectors require progressive precision, and thus have superlinear storage cost.
  This section addresses how \cref{alg:cfas-naive,alg:fmg-residual-naive} can be implemented instead with $\mathcal O(n)$ memory.

  \subsection{Matrices and \texorpdfstring{\gls*{rhs}}{RHS}}

  The storage cost of matrices becomes irrelevant if the method is implemented matrix-free.
  Whenever memory footprint is a concern, a matrix-free implementation should be considered before thinking about reducing the storage cost of vectors (the scope of this work), as matrices typically require significantly more memory than vectors.
  Moreover, matrix-free methods have been shown to outperform matrix-based approaches due to the typical imbalance of memory bandwidth and compute performance of modern hardware~\cite{kronbichler_performance_2018}.

  Likewise, an analytically given \gls*{rhs} can be assembled on-the-fly so that it requires practically no memory.
  If the \gls*{rhs} data is not known analytically but instead stored on disk, it can be read in streaming fashion, without ever loading it completely to memory.
  In case the \gls*{rhs} arises from the solution of a different \gls*{pde} and overall linear space complexity is desired, the \gls*{rhs} must be computed and stored in compact format as well.
  A compact \gls*{rhs} can be decoded in the coarse-to-fine sweep of the residual calculation (\cref{alg:fmg-residual-naive}) exactly like the solution.
  All mentioned approaches lead to linear memory footprint.

  \begin{remark}
    The scope of this work is developing a solver capable of solving a \gls*{pde} to arbitrary accuracy.
    In this scenario, the linear operator and \gls*{rhs} must inevitably be known to arbitrary accuracy, e.g., from an analytical description.
    Therefore, we focus herein on the case of an analytical \gls*{rhs}, and in particular choose this setting for the numerical experiments below.
    This allows for on-the-fly assembly of the \gls*{rhs} and eliminates the need of storing the full \gls*{rhs} vector.
    In practice, the solution is only interesting to a certain accuracy, which limits the required accuracy of the operator and \gls*{rhs}.
    In particular, a typical engineering / domain science application uses coefficient and \gls*{rhs} data that is not exactly known.
    Common sources of uncertainty are modelling and measurement error.
    It is important to realize that a solver that does not balance all error sources is inefficient.
    Solving the \gls*{lse} more accurately does not result in a more accurate solution to the \gls*{pde} if the discretization error is large.
    Similarly, solving the \gls*{pde} more accurately by computing on a finer mesh does not result in a more accurate solution to the engineering problem if the \gls*{rhs} is a crude approximation of reality.
    Therefore, it is crucial in practice to estimate the magnitude of each of these terms and be aware of the target solution accuracy; otherwise, an adequate grid spacing cannot be chosen.
    Ultimately, computing on a needlessly large grid where the solution accuracy is overall limited by errors in the \gls*{pde} (modelled physics, coefficient functions, \gls*{rhs}) would be a waste of electric energy.
  \end{remark}

  \subsection{Temporary vectors}

  The key idea behind matrix-free methods is to never store the entire matrix but only a few select values that are required to determine the result of the operator application at the current \gls*{dof}/element.
  Similarly, we limit the memory consumption of the temporary vectors $\std z_{0\rng L}$, $\std u_{0\rng L}$, and $\std t_{0\rng L}$ by minimizing the number of entries that are simultaneously present in memory.

  Consider the coarse-to-fine sweep of \cref{alg:cfas-naive} as an example.
  It iterates over the grid hierarchy computing $\std z_{0\rng L}$ and from it $\compact y_{0\rng L}$.
  Crucially, for all levels $\ell$, $\std z_\ell$ is needed only temporarily, to compute $\compact y_\ell$ and $\std z_{\ell+1}$.
  As soon as these dependent quantities have been derived, $\std z_\ell$ can be evicted from memory.
  The time each individual entry of $\std z_\ell$ must reside in memory, and thus the total number of values stored at the same time, can be minimized by traversing the grid hierarchy in a specific order.

  To find this ordering, we exploit the locality of stencil applications, i.e., the sparsity pattern of the matrices.
  Due to locality, the computation of $\compact y_\ell$ and $\std z_{\ell+1}$ can begin before all entries of $\std z_\ell$ have been computed.
  The key strategy is to compute dependent values as soon as possible with the available data in $\std z_\ell$.
  This loosely corresponds to a \gls*{dfs} traversal of the dependency graph, while the naive variant is closely related to \gls*{bfs}.
  Also due to locality, a certain \gls*{dof} in $\std z_\ell$ only influences a few values in the dependent vectors.
  Hence, this \gls*{dof} can be evicted from memory as soon as all resulting values have been derived from it.
  Consequently, only a small portion of $\std z_\ell$ must reside in memory at a given time.

  As a result of the above, $\std z_{0\rng L}$ and $\compact y_{0\rng L}$ are computed in a local region on \emph{all} grids, instead of on the entire domain on one specific level.
  Owing to the locality of stencil applications, this requires only data from a small neighborhood.
  Moreover, only the boundary of the already handled area need be reused in future iterations and must therefore be stored.
  Because the boundary has codimension one, it covers only $\mathcal O(n_\ell^{(d-1)/d})$ points.

  More precisely, the number of points that must be stored to compute a matrix-vector product $M \std x$, $M \in \mathbb R^{m\times n}$, depends on the matrix' sparsity pattern.
  Denote by $\fstCol(M,\rho)$ and $\lstCol(M,\rho)$ the minimum and maximum column index of the non-zero entries in row $\rho$, respectively.
  To apply the stencil of $M$ at grid point~$\rho$, only \glspl*{dof} of $\std x$ with indices between $\fstCol(M,\rho)$ and $\lstCol(M,\rho)$ are required.
  To avoid recomputing \glspl*{dof} required in future iterations at grid points $\sigma > \rho$, only \glspl*{dof} with indices smaller than $\min_{\rho \leq \sigma \leq m} \fstCol(M,\sigma)$ are evicted from memory.
  Therefore, at most
  \begin{equation}
    \label{eq:tail-len}
    \tail(M \in \mathbb R^{m\times n}) := \max_{1 \leq \rho \leq m} ( \lstCol(M,\rho) - \min_{\rho \leq \sigma \leq m} \fstCol(M,\sigma) + 1 )
  \end{equation}
  entries of the vector $\std x$ must be stored at a time.
  Clearly, $\tail$ depends on the ordering of the unknowns.
  On structured grids with lexicographic ordering,
  \begin{equation}
    \label{eq:tail-len-lexicographic}
    \tail(M_\ell) \approx (w(M_\ell)-1) n_\ell^{\frac{d-1}{d}}\,,
  \end{equation}
  where $w(M_\ell)$ is the stencil width of $M_\ell$, e.g., 3 for a 3-point stencil in 1D, a 5-point stencil in 2D, or a 7-point stencil in 3D.
  We show in \cref{sec:cost-temporaries} that the $\delta$ in \cref{eq:tail-len-lexicographic} is sufficiently small to guarantee linear space complexity.
  In any case, a clever ordering of unknowns can be used to minimize $\tail$.
  For example, orderings based on space-filling curves come to mind.
  However, one has to be careful since these typically do not minimize for the \emph{maximum} distance of neighboring nodes, which \cref{eq:tail-len} is based on.
  In that case, it would be best to limit the cache to values that will actually be reused, instead of storing everything with an index greater or equal to the smallest reused index.
  Alternatively, a trade-off between buffer size and number of recomputations can be made.

  \begin{algorithm}[tbp]%
    \algosize%
    \subfloat[First, the loop (\cref{line:cfas-c2f-loop} in \cref{alg:cfas-naive}) is transformed to an equivalent recursive form (highlighted lines).]{%
      \label{alg:c2f-rec}%
      \begin{pseudo}[kw,pos=b]%
        $\std z_0 \gets 0$ \\
        $\compact y_0 \gets \compact y_0 + M_0 (\std r_0 - \std s_0 - A_0 \compact y_0)$ \\[hl]
        if $L > 0$ \\+[hl]
          $(\std z, \compact y)_{1\rng L} \gets \pr{c2f}(\std z_0, \compact y_{0\rng L}, 1)$ \\-*[hl]

        \addlinespace[2ex]
        \midrule
        \addlinespace[2ex]

        \pseudohd{c2f}(\std z_{\ell-1}, \compact y_{\ell-1 \rng L}, \ell) \\
        $\std z_\ell \gets P_\ell \compact y_{\ell-1} + P_\ell \std z_{\ell-1}$ \label{line:c2f-rec-z} \\
        $\compact y_\ell \gets \compact y_\ell + M_\ell (\std r_\ell - \std s_\ell - A_\ell \std z_\ell - A_\ell \compact y_\ell)$\label{line:c2f-rec-v}\\[hl]

        if $\ell < L$ \\+[hl]
          $(\std z, \compact y)_{\ell+1 \rng L} \gets \pr{c2f}(\std z_\ell, \compact y_{\ell\rng L}, \ell+1)$ \\-[hl]
        return $(\std z, \compact y)_{\ell\rng L}$ \label{line:c2f-rec-end}
      \end{pseudo}%
    }\hfill
    \subfloat[Next, the vector operations (\cref{line:c2f-rec-z,line:c2f-rec-v} in \cref{alg:c2f-rec}) are written as loops over the grid wherein a stencil application is performed at each node.]{%
      \label{alg:c2f-rec-loops}%
      \begin{pseudo}[kw,pos=b]%
        $\std z_0 \gets 0$ \\
        $\compact y_0 \gets \compact y_0 + M_0 (\std r_0 - \std s_0 - A_0 \compact y_0)$ \\[hl]
        $i_0, j_0 \gets n_0$ \label{line:c2f-rec-loops:ij0} \\[hl]
        $i_{1\rng L}, j_{1\rng L} \gets 0_{1\rng L}$ \label{line:c2f-rec-loops:ijl} \\
        if $L > 0$ \\+
          $(\std z, i, \compact y, j)_{1\rng L} \gets \pr{c2f}(\std z_0, i_{1\rng L}, \compact y_{0\rng L}, j_{1\rng L}, 1)$ \\-*

        \addlinespace[2ex]
        \midrule
        \addlinespace[2ex]

        \pseudohd{c2f}(\std z_{\ell-1}, i_{\ell\rng L}, \compact y_{\ell-1 \rng L}, j_{\ell\rng L}, \ell) \\[hl]
        while $i_\ell < n_\ell$ \label{line:c2f-rec-loops:i-loop} \\+
          $\std z_{\ell,i_\ell} \gets (P_\ell \compact y_{\ell-1})_{i_\ell} + (P_\ell \std z_{\ell-1})_{i_\ell}$ \\[hl]
          $i_\ell \gets i_\ell+1$ \\-[hl]

        while $j_\ell < n_\ell$ \label{line:c2f-rec-loops:j-loop} \\+
          $\compact y_{\ell,j_\ell} \gets (\compact y_\ell + M_\ell (\std r_\ell - \std s_\ell - A_\ell \std z_\ell - A_\ell \compact y_\ell))_{j_\ell}$ \\[hl]
          $j_\ell \gets j_\ell + 1$ \label{line:c2f-rec-loops:j-loop-end} \\-

        if $\ell < L$ \\+
          $(\std z, i, \compact y, j)_{\ell+1 \rng L} \gets \pr{c2f}(\std z_\ell, i_{\ell+1 \rng L}, \compact y_{\ell\rng L}, j_{\ell+1 \rng L}, \ell+1)$ \label{line:c2f-rec-loops:rec} \\-
        return $(\std z, i, \compact y, j)_{\ell\rng L}$
      \end{pseudo}%
    }

    \caption{%
      Intermediate steps for transforming the naive coarse-to-fine sweep (\crefrange{line:cfas-n-z0}{line:cfas-n-yl}) of \cref{alg:cfas-naive} to the optimized version in \cref{alg:cfas-opt-recursive}.%
    }
  \end{algorithm}

  \begin{algorithm}[tp]%
    \algosize%
    \caption{\pr{cFas}: \Cref{alg:cfas-naive} optimized for linear space complexity.}%
    \label{alg:cfas-opt-recursive}%
%
    \addtocounter{algorithm}{-1}%
    \paragraph{Input}
    \begin{tabular}{ll}
      \precA{$A_{0\rng L}$}      & matrix hierarchy \\
      \precA{$M_{0\rng L}$}      & smoother \\
      \precA{$P_{1\rng L}$}      & prolongation matrices \\
      \precR{$\compact y_{0\rng L} = 0_{0\rng L}$} & initial guess of solution \\
      \precR{$\std r_{0\rng L}$} & \gls*{rhs} \\
    \end{tabular}

    \paragraph{Output}
    Solution $\precR{\compact y_{0\rng L}}$ improved by one V(0,1) multigrid cycle.
    \unskip\medskip

    \subfloat[Main algorithm.]{%
      \begin{pseudo}[fullwidth,kw]*&
        \ct{Fine-to-coarse sweep is omitted because \precR{$\compact y_{0\rng L} = 0_{0\rng L}$}.} \\*&
        \ct{Coarse-to-fine sweep.} \\

        $\precZ{\std z_0 \gets 0}$ \\
        $\precR{\compact y_0 \gets \flFix{\precM{\flMaxFulPRO{\precR{\compact y_0}} + \flMaxFulPRO{\precA{M_0}} (\flMaxFulPRO{\precR{\std r_0}} - \flMaxFulPRO{\precA{A_0}} \flMaxFulPRO{\precR{\compact y_0}})}}}$ \hfill\ct{Smooth $A_0 \compact y_0 = \std r_0$} \\

        $i_0, j_0 \gets n_0$ \hfill\ct{$i$: index of next node in $\std z$.} \\
        $i_{1\rng L}, j_{1\rng L} \gets 0_{1\rng L}$ \hfill\ct{$j$: index of next node in $\std v$.} \\

        for $\ell = 1..(L-1)$ step $1$ \label{line:cfas-opt-rec:init-z-start} \\+
        $\precZ{\std z_\ell} \gets{}$ \tn{Empty buffer with a capacity of $\max\{\tail(A_\ell),\tail(P_{\ell+1}),\tail(A_{\ell+1}P_{\ell+1})\}$ \glspl*{dof}.} \\-
        $\precZ{\std z_L} \gets{}$ \tn{Empty buffer with a capacity of $\tail(A_L)$ \glspl*{dof}.} \label{line:cfas-opt-rec:init-z-end} \\

        if $L > 0$ \\+
          $(\precZ{\std z}, i, \precR{\compact y}, j)_{1\rng L} \gets \pr{c2f}({(\precZ{\std z}, i, \precR{\compact y}, j)_{0\rng L}}, 1)$ \\-

        return \precR{$\compact y_{0\rng L}$}
      \end{pseudo}%
    }%
    \unskip\medskip

    \subfloat[\pr{c2f}({(\precZ{\std z}, i, \precR{\compact y}, j)_{\ell-1 \rng L}}, \ell): Recursive coarse-to-fine sweep.]{%
      \begin{pseudo}[fullwidth,kw,start=11]
        while $i_\ell < n_\ell$ and $\lstCol(P_\ell, i_\ell) < j_{\ell-1}$ and $\lstCol(P_\ell, i_\ell) < i_{\ell-1}$ \hfill\ct{For all nodes $i_\ell$ where $\std z_\ell$ can be computed.} \label{line:cfas-rec-c2f:i-loop} \\+

          if \tn{capacity of $\precZ{\std z_\ell}$ is exhausted} drop \tn{oldest value of $\precZ{\std z_\ell}$} \label{line:cfas-rec-c2f:drop-z} \hfill\ct{Limit memory consumption.} \\

          $\precZ{\std z_{\ell,i_\ell} \gets \flPro{\precM{\precM{(\flMaxFulPRO{\precA{P_\ell}} \flMaxFulPRO{\precR{\compact y_{\ell-1}}})_{i_\ell} + (\flMaxFulPRO{\precA{P_\ell}} \flMaxFulPRO{\precZ{\std z_{\ell-1}}})_{i_\ell}}}}}$ \label{line:cfas-rec-c2f:apply-P} \\
          $i_\ell \gets i_\ell+1$ \\

          while $j_\ell < n_\ell$ and $\lstCol(A_\ell, j_\ell) < i_\ell$ \label{line:cfas-rec-c2f:j-loop} \hfill\ct{For all nodes $j_\ell$ where $\compact y_\ell$ can be computed.} \\+
            $\precR{\compact y_{\ell,j_\ell} \gets \flFix{\precM{(\flMaxFulPRO{\precR{\compact y_\ell}} + \flMaxFulPRO{\precA{M_\ell}} (\flMaxFulPRO{\precR{\std r_\ell}} - \flMaxFulPRO{\precA{A_\ell}} \flMaxFulPRO{\precZ{\std z_\ell}} - \flMaxFulPRO{\precA{A_\ell}} \flMaxFulPRO{\precR{\compact y_\ell}}))_{j_\ell}}}}$ \hfill\ct{Smooth $A_\ell \compact y_\ell = \std r_\ell - A_\ell \std z_\ell$} \\
            $j_\ell \gets j_\ell + 1$ \label{line:cfas-rec-c2f:incr-j} \\-

          if $\ell < L$ \\+*&
            \ct{Process all nodes on finer grids that can be computed with the available data.} \\
            $(\precZ{\std z}, i, \precR{\compact y}, j)_{\ell+1 \rng L} \gets \pr{c2f}({(\precZ{\std z}, i, \precR{\compact y}, j)_{\ell\rng L}}, \ell+1)$ \label{line:cfas-rec-c2f:rec} \\--

        return $(\precZ{\std z}, i, \precR{\compact y}, j)_{\ell\rng L}$
      \end{pseudo}%
    }%
    \stepcounter{algorithm}%
  \end{algorithm}

  \begin{algorithm}[pt]%
    \algosize%
    \caption{\pr{fmgResidual}: \Cref{alg:fmg-residual-naive} optimized for linear space complexity.}%
    \label{alg:fmg-residual-opt-recursive}%
%
    \addtocounter{algorithm}{-1}%
    \paragraph{Input}
    \begin{tabular}{ll}
      \precF{$A_L$}                   & matrix \\
      \precF{$P_{1\rng L}$}           & prolongation matrices \\
      \precF{$R_{1\rng L}$}           & restriction matrices \\
      \precF{$\std f_L$}              & \gls*{rhs} \\
      \precU{$\compact u_{0 \rng L}$} & solution approximation \\
    \end{tabular}

    \paragraph{Output}
    Residual $\precR{\std r_\ell} = R_{\ell+1} \cdots R_L (\std f_L - A_L (\compact u_L + \cdots + P_L \cdots P_1 \compact u_0))$ for $0 \leq \ell \leq L$.
    \unskip\medskip

    \subfloat[Main algorithm.]{%
      \begin{pseudo}[fullwidth,kw]
        $\std u_0 \gets \precU{\compact u_0}$ \\
        $i_0 \gets n_0$ \hfill\ct{$i$: index of next node in $\std u$.} \\
        $i_{1\rng L} \gets 0_{1\rng L}$ \\
        $j_{0\rng L} \gets 0_{0\rng L}$ \hfill\ct{$j$: index of next node in $\std t$.} \\

        for $\ell = 1..(L-1)$ step $1$ \\+
        $\std u_\ell \gets{}$ \tn{Empty buffer with a capacity of $\tail(P_{\ell+1})$ \glspl*{dof}.} \\-
        $\std u_L \gets{}$ \tn{Empty buffer with a capacity of $\tail(A_L)$ \glspl*{dof}.} \\

        $\std t_0 \gets{}$ \tn{Empty buffer with a capacity of $1$ \gls*{dof}.} \\
        for $\ell = 1..L$ step $1$ \\+
        $\std t_\ell \gets{}$ \tn{Empty buffer with a capacity of $\tail(R_\ell)$ \glspl*{dof}.} \\-

        if $L > 0$ \\+
          $(\std u, i)_{1\rng L}, (\std t, j, \precR{\std r})_{0\rng L} \gets \pr{c2f}(\precU{\compact u_{1 \rng L}}, {(\std u, i)_{0 \rng L}}, {(\std t, j, \precR{\std r})_{0\rng L}}, 1)$ \\-
        return \precR{$\std r_{0 \rng L}$}
      \end{pseudo}%
    }%
    \unskip\medskip

    \subfloat[\pr{c2f}(\precU{\compact u_{\ell \rng L}}, {(\std u, i)_{\ell-1 \rng L}}, {(\std t, j, \precR{\std r})_{0\rng L}}, \ell): Recursive coarse-to-fine sweep.]{%
      \begin{pseudo}[fullwidth,kw,start=14]
        while $i_\ell < n_\ell$ and $\lstCol(P_\ell, i_\ell) < i_{\ell-1}$ \hfill\ct{For all nodes $i_\ell$ where $\std u_\ell$ can be computed.} \\+

          if \tn{capacity of $\std u_\ell$ is exhausted} drop \tn{oldest value of $\std u_\ell$} \hfill\ct{Limit memory consumption.} \\

          $\std u_{\ell,i_\ell} \gets \flFul{\precU{\compact u_{\ell,i_\ell}}} + \flFul{\precM{(\flMaxFulPRO{\precF{P_\ell}} \flMaxFulPRO{\precL{\std u_{\ell-1}}})_{i_\ell}}}$ \\
          $i_\ell \gets i_\ell+1$ \\

          if $\ell = L$ \\+
            $(\std t, j, \precR{\std r})_{0\rng L} \gets \pr{f2c}({(\std t, j, \precR{\std r})_{0\rng L}}, L)$ \\-

          else \\+*&
            \ct{Process all nodes on finer grids that can be computed with the available data.} \\
            $(\std u, i)_{\ell+1 \rng L}, (\std t, j, \precR{\std r})_{0\rng L} \gets \pr{c2f}(\precU{\compact u_{\ell+1 \rng L}}, {(\std u, i)_{\ell \rng L}}, {(\std t, j, \precR{\std r})_{0\rng L}}, \ell+1)$ \\--

        return $(\std u, i)_{\ell\rng L}, (\std t, j, \precR{\std r})_{0\rng L}$
      \end{pseudo}%
    }%
    \unskip\medskip

    \subfloat[\pr{f2c}({(\std t, j, \precR{\std r})_{0\rng \ell+1}}, \ell): Recursive fine-to-coarse sweep.]{%
      \begin{pseudo}[fullwidth,kw,start=23]
        $\id{canApplyStencil} \gets \begin{cases}
          \lstCol(A_L, j_L) < i_L                  & \text{if } \ell = L \\
          \lstCol(R_{\ell+1}, j_\ell) < j_{\ell+1} & \text{else}
        \end{cases}$ \\

        while $j_\ell < n_\ell$ and \id{canApplyStencil} \\+

          if \tn{capacity of $\std t_\ell$ is exhausted} drop \tn{oldest value of $\std t_\ell$} \\

          $\std t_{\ell,j_\ell} \gets \begin{cases}
            \flFul{\precM{\flMaxFulPRO{\precF{\std f_{L,j_L}}} - (\flMaxFulPRO{\precF{A_L}} \flMaxFulPRO{\precL{\std u_L}})_{j_L}}} & \text{if } \ell = L \\
            \flFul{\precM{(\flMaxFulPRO{\precF{R_{\ell+1}}} \flMaxFulPRO{\precL{\std t_{\ell+1}}})_{j_\ell}}} & \text{else}
          \end{cases}$ \\
          $\precR{\std r_{\ell,j_\ell} \gets \flFix{\precL{\std t_{\ell,j_\ell}}}}$ \\
          $j_\ell \gets j_\ell + 1$ \\

          if $\ell > 0$ \\+
            $(\std t, j, \precR{\std r})_{0\rng \ell-1} \gets \pr{f2c}({(\std t, j, \precR{\std r})_{0\rng \ell}}, \ell-1)$ \\--

        return $(\std t, j, \precR{\std r})_{0\rng \ell}$
      \end{pseudo}%
    }%
    \stepcounter{algorithm}%
  \end{algorithm}

  We now apply the above discussion and derive the memory-optimal implementation of the coarse-to-fine sweep in \cref{alg:cfas-naive} through a series of code transformations.
  First, the iterative sweep can be equivalently written recursively as in \cref{alg:c2f-rec}.
  For any level index $\ell$, $0 < \ell \leq L$, the recursive \pr{c2f} routine (\crefrange{line:c2f-rec-z}{line:c2f-rec-end}) computes and returns $\std z_{\ell\rng L}$ and $\compact y_{\ell\rng L}$.
  Level $\ell$ is handled directly, and the recursion fills in the finer grids.
  Analogous to the loop in \cref{line:cfas-c2f-loop} of \cref{alg:cfas-naive}, the index $\ell$ is incremented each recursive invocation until the base case $\ell = L$ is reached.

  The next transformation is to write the vector operations in \cref{line:c2f-rec-z,line:c2f-rec-v} of \cref{alg:c2f-rec} as explicit loops over the grid (\cref{alg:c2f-rec-loops}, \crefrange{line:c2f-rec-loops:i-loop}{line:c2f-rec-loops:j-loop-end}).
  To that end, arrays of indices $i_{0\rng L}$ and $j_{0\rng L}$ are introduced (\cref{line:c2f-rec-loops:ij0,line:c2f-rec-loops:ijl}).
  The index $i_\ell$, with $0 \leq i_\ell \leq n_\ell$, points to the next \gls*{dof} in $\std z_\ell$ that must be computed.
  Likewise, $j_\ell$ counts the progress of the computation of $\compact y_\ell$.
  Through the loops in \crefrange{line:c2f-rec-loops:i-loop}{line:c2f-rec-loops:j-loop-end}, the respective stencil is applied at each node and the index incremented accordingly.
  We denote the stencil application of, e.g., $P_\ell$ to $\compact y_{\ell-1}$ at node $i_\ell$ as $(P_\ell \compact y_{\ell-1})_{i_\ell}$.

  The crucial step allowing for the memory-optimization is to move the computation of $\compact y_\ell$ and the recursive call (\crefrange{line:c2f-rec-loops:j-loop}{line:c2f-rec-loops:rec} of \cref{alg:c2f-rec-loops}) into the $i_\ell$-loop, as shown in the optimized version \cref{alg:cfas-opt-recursive}, \crefrange{line:cfas-rec-c2f:j-loop}{line:cfas-rec-c2f:rec}.
  This way, as soon as a \gls*{dof} of $\std z_\ell$ is calculated (\cref{line:cfas-rec-c2f:apply-P} of \cref{alg:cfas-opt-recursive}), the computation directly proceeds with the dependent values of $\compact y_\ell$ and the finer grids.
  Only then does the iteration continue with the next node of $\std z_\ell$.
  Because only certain entries of the dependent vectors can be computed, we introduce additional bound checks in \cref{line:cfas-rec-c2f:i-loop,line:cfas-rec-c2f:j-loop} of \cref{alg:cfas-opt-recursive}.
  Finally, this change in the order in which the nodes are processed allows \cref{line:cfas-rec-c2f:drop-z} to evict \glspl*{dof} of $\std z_\ell$ quickly from memory.
  Therefore, in \crefrange{line:cfas-opt-rec:init-z-start}{line:cfas-opt-rec:init-z-end}, we initialize buffers for $\std z_{1\rng L}$ with the minimally possible capacity such that the stencils of $A$ and $P$ can be applied at each grid point.
  As noted above, we assume that the initial guess for $\compact y$ is zero, so that $\std s = 0$ and no work has to be done in the fine-to-coarse sweep of \cref{alg:cfas-opt-recursive}.

  \Cref{alg:fmg-residual-naive} is optimized in a very similar fashion.
  The optimized version is displayed in \cref{alg:fmg-residual-opt-recursive}.
  The key difference to \cref{alg:cfas-naive} is that the coarse-to-fine sweep is followed by a fine-to-coarse sweep.
  This manifests in the optimized implementation (\cref{alg:fmg-residual-opt-recursive}) in that the fine-to-coarse sweep is nested in the coarse-to-fine sweep.
  Apart from that, the pattern remains the same.
  Note that we presented recursive forms of the optimized algorithms for ease of derivation, but they have equivalent iterative versions.

  \section{Cost analysis}
  \label{sec:cost}

  In this section, we discuss the cost associated with the above algorithms.
  Storage, computation, and memory transfers are considered.

  \subsection{Storage}

  We predicate that the implementation is matrix-free so that the cost to store stencils/matrices and the \gls*{rhs} can be ignored.
  Below, we show that mantissas and exponents of all vectors can be stored with $\mathcal O(n_L)$ bits.

  \subsubsection{Mantissas of regressive-precision vectors}

  Starting with the re\-gres\-sive-precision vectors $\compact u_{0\rng L}$, $\std r_{0\rng L}$, and $\compact y_{0\rng L}$, assume that the number of \glspl*{dof} on level~$\ell$ is $n_\ell \approx c 2^{d\ell}$, with a positive constant $c$ possibly dependent of $p$.
  Recall that a $\regressive{k}{b}$ precision vector comprises $k (L-\ell) + b$ bits on level~$\ell$.
  Multiplying by $n_\ell$ and summing over all grids, we find the total number of bits used:\pagebreak[0]
  \begin{equation}
    N_d^{\compact u, \std r, \compact y}(k,b) \approx c \sum_{\ell=0}^L (k(L-\ell)+b) 2^{d \ell}
    \overset{d\neq0}{<} \left( \frac{2^d}{2^d-1} b + \frac{2^d}{(2^d-1)^2} k \right) n_L \in \mathcal O(n_L)\,.
  \end{equation}
  This establishes linear storage cost w.r.t.\ the number of unknowns.
  The constant $2^d (2^d-1)^{-1}$ is $2$ in 1D, $\frac 43$ in 2D, and $\frac 87$ in 3D, and is therefore small.

  \subsubsection{Mantissas of temporary vectors}
  \label{sec:cost-temporaries}

  The remaining vectors $\std z_{0\rng L}$, $\std u_{0\rng L}$, and $\std t_{0\rng L}$ require progressive precision, but only $\delta (M)$ (cf.\ \cref{eq:tail-len}) entries are stored in the optimized \cref{alg:cfas-opt-recursive,alg:fmg-residual-opt-recursive}.
  Here, $M$ is the matrix that maximizes $\delta(M)$ among those applied to the respective temporary vector.
  Concretely, this is $A$ for $\std u$ and $\std z$ and $R$ for $\std t$ when using B-splines for the finite element discretization.
  Assuming the unknowns constitute a uniform grid and are ordered lexicographically, it has been established in \cref{sec:memory-optimization} that $\delta (M_\ell) \approx (w(M_\ell) - 1) n_\ell^{(d-1)/d}$.

  Further assuming that the stencil width of the relevant matrices is independent of $\ell$, the storage cost of $\std u$ and $\std t$ in \cref{alg:fmg-residual-opt-recursive} and $\std z$ in \cref{alg:cfas-opt-recursive} are given by \cref{eq:cost-tmp-fix} and \cref{eq:cost-tmp-pro}, respectively, where $a_L := kL+b$ and $\tilde c_\ell := (w(M_\ell) - 1)$.
  In all cases, the space complexity is strictly less than $\mathcal O(n_L)$.\pagebreak[0]
  \begin{gather}
    \begin{aligned}
      \label{eq:cost-tmp-fix}
      \std N_d^{\std u, \std t}(a_L)
      &= \sum_{\ell=0}^L a_L \tail(M_\ell)
      \approx \tilde c a_L \sum_{\ell=0}^L n_\ell^{\frac{d-1}{d}} \\
      &\leq \begin{cases}
        \tilde c a_L (L+1) \in \mathcal O(\log^2 n_L), & d = 1 \\
        \frac{2^{d-1}}{2^{d-1}-1} \tilde c a_L n_L^{\frac{d-1}{d}} \in \mathcal O(n_L^{\frac{d-1}{d}} \log n_L), & d > 1
      \end{cases}
    \end{aligned} \displaybreak[0]\\
    \begin{aligned}
      \label{eq:cost-tmp-pro}
      \std N_d^{\std z}(k,b)
      &= \sum_{\ell=0}^L (k\ell+b)\tail(M_\ell)
      \leq \sum_{\ell=0}^L a_L \tail(M_\ell)
      = \std N_d^{\std u, \std t}(a_L)
    \end{aligned}
  \end{gather}

  \subsubsection{Exponents}

  Apart from the mantissa, a floating point number also carries an exponent.
  The exponents of the compact solution vector decrease by $p+1$ for each refinement level.
  Consequently, storing one exponent for each \gls*{dof} would require $\mathcal O(n_L \log \log n_L)$ memory.

  Therefore, we instead store $\compact u_\ell$, $\std r_\ell$, and $\compact y_\ell$ for all $0 \leq \ell \leq L$ in \gls*{bfp} format~\cite{wilkinson_rounding_2023}, rendering the storage cost of exponents negligible.
  In \gls*{bfp} format, the vector is divided into few blocks.
  Each block comprises several mantissas, but only a single exponent scaling all mantissas.
  The extreme case, in which only a single exponent is stored for the entire vector, is also known as scaled fixed-point format.
  The applicability of \gls*{bfp} to mixed-precision multigrid methods has been demonstrated in~\cite{kohl_multigrid_2024}.

  One of the difficulties of this format is that the maximal exponent must be known before the mantissas can be normalized.
  In practice, the exponent might be known/estimated a priori~\cite{kohl_multigrid_2024}.
  Depending on how good the estimate is, the mantissa's size can be temporarily increased to allow for an inaccurate guess.
  This technique typically requires normalization after the computation, effectively doubling the memory traffic.
  If no (good) estimation is available, the values can be computed twice:
  once to determine the maximum exponent and a second time to store the values (according to the \gls*{bfp} exponent)~\cite{kohl_multigrid_2024}.

  Both re-normalization and re-computation can be avoided by increasing the number of blocks/exponents to $b$, with $b$ being the number of bits in each mantissa.
  As the entries of the vector are computed, we keep track of the currently largest known exponent and normalize the current mantissa to that exponent.
  This exponent and associated mantissas constitute one block of the \gls*{bfp} vector.
  The next block starts when a larger exponent is encountered.
  Finally, to make the storage cost of the exponents independent of the vector size, note that storing only the $b$ largest exponents suffices to recover all entries to single-exponent \gls*{bfp} accuracy:
  The other exponents are smaller by at least $b$ compared to the largest exponent.
  Therefore, those entries would be truncated to zero in single-exponent \gls*{bfp} precision, irrespective of the actual value of the exponent or mantissa.
  Hence, at most $b$ exponents and block-sizes must be stored to match the accuracy of the scaled fixed-point format.

  To achieve linear space complexity, we use the \gls*{bfp} format for $\compact u_\ell$, $\std r_\ell$, and $\compact y_\ell$.
  Stencil entries and the temporary vectors $\std z_\ell$, $\std u_\ell$, and $\std t_\ell$ have linear space complexity even in \gls*{fp} format because only few entries reside in memory at a time.
  Concretely, we found that the temporaries' exponents also grow linearly with $L$, and can therefore be stored with $\mathcal O(n_L^{(d-1)/d} \log \log n_L) \subsetneq \mathcal O(n_L)$ bits.
  To avoid overcomplicating matters, we use the usual \gls*{fp} format for stencils and temporaries.
  Therefore, whenever one of the arguments to an arithmetic operation is stored in \gls*{fp} format, the \gls*{bfp} argument is converted to its \gls*{fp} representation before performing the operation in \gls*{fp} arithmetic.
  In particular, all stencil applications use \gls*{fp} arithmetic and entail the associated rounding behavior.
  While we use \gls*{bfp} only where it is required for linear space complexity, storing and computing everything in \gls*{bfp} format as in~\cite{kohl_multigrid_2024} has the additional advantage of relying solely on integer arithmetic, which is cheaper than \gls*{fp} arithmetic.

  \subsubsection{Comparison to progressive precision \texorpdfstring{\gls*{fmg}}{FMG} in \texorpdfstring{\gls*{bfp}}{BFP} arithmetic}

  To demonstrate that the improved asymptotics pay off at realistic problem sizes, we estimate the memory consumption of the isolated solution vector and an entire \gls*{fmg} implementation using the formulas from above.
  We model the storage cost on paper instead of measuring the memory footprint on an actual machine because, for flexibility’s sake, our implementation emulates all low precision numbers using arbitrary precision libraries optimized for large bit widths, and each value consumes at least 416 bits.

  The baseline for the comparison is \gls*{fmg} in progressive precision \gls*{bfp} arithmetic (BFPMG)~\cite{kohl_multigrid_2024}.
  To the best of our knowledge, this is the current state-of-the-art in terms of memory efficiency.
  We make the following assumptions:
  BFPMG is implemented in a memory-efficient manner.
  In particular, the coarse-grid solution vector is freed immediately after computing the \gls*{fmg} prolongation.
  Moreover, during the V-cycle, memory is required only for the vectors $\std u_L$, $\std r_L$, $\std r_{L-1}$, and $\std y_{0\rng L}$.
  As in \gls*{cmg}, matrices and \gls*{rhs} are assembled on-the-fly and do not consume memory.
  When the matrix-free method is efficiently implemented, except for the vectors, no additional $h$-dependent data structures need to be stored~\cite{kohl_fundamental_2023}.
  Lastly, we assume that no \gls*{bfp} normalizations are performed within the \gls*{gemv} routines and thus no additional storage is required in \gls*{gemv} operations nor the smoother.

  Obviously, the memory consumed by either method is determined by the precision used to represent the vectors.
  While the precision \emph{increments} are directly determined by $m$ and $p$, the precision \emph{constants} depend in a less predictable way on specifics like the domain and \gls*{rhs}.
  Instead of guessing, we determine these constants for some actual model problems: Poisson's equation in 1D and 2D ($m=1$) and the biharmonic equation in 1D ($m=2$).
  Details about these \glspl*{pde} and how we determine the precision constants are given in \cref{sec:experiments}.
  The constants are listed in \cref{tab:experiments-bits}.
  This gives us all we need to precisely predict the required memory of \gls*{cmg}.
  For BFPMG, we simply assume that the same constants apply, which is why we refer to the saving figures as estimates.
  Note, though, that in BFPMG, all vectors are stored in progressive precision, meaning that storage requirements are ultimately dominated by the precision increments while the significance of the constants diminishes with increasing number of levels.

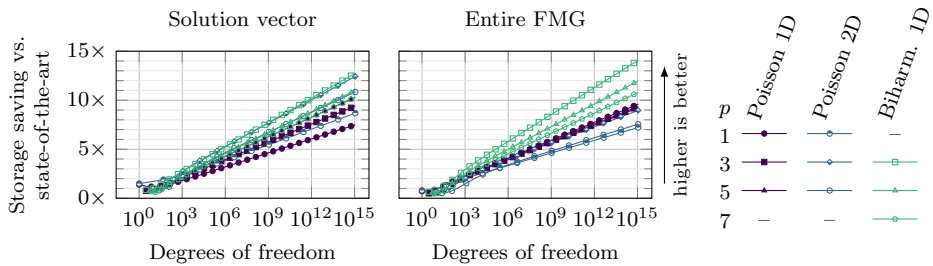
\begin{figure}
    \centering%
    \subfloat{%
    \begin{tikzpicture}%
      \begin{groupplot}[
        group style = {group size=3 by 2, x descriptions at=edge bottom, y descriptions at=edge left, horizontal sep=4mm, vertical sep=6mm},
        width = 0.36\linewidth,
        xlabel = {\Acrlongpl*{dof}},
        every axis/.append style = {
          xmode = log,
          xtickten = {0,3,...,15},
          minor y tick num = 4,
          yticklabel = {$\pgfmathprintnumber{\tick}\times$},
          xmax = 10^16,
          ymin = 0,
          ymax = 15,
        },
        cycle multiindex list = {[samples of colormap = 8 of viridis]\nextlist mymarklist},
        /tikz/mark repeat = 2,
      ]
        \nextgroupplot[
          title = {Poisson 1D},
          ylabel = {Entire \gls*{fmg}\\[1.5ex]Storage saving vs.\\state-of-the-art},
          legend to name = {legend:savings},
          legend columns = 8,
        ]
        \addplot table[x=dofs-m1-d1-p01,y=prog-m1-d1-p01/cReg] {data/memory-saving-fmg.dat}; \addlegendentry{$p=1$}
        \addplot table[x=dofs-m1-d1-p02,y=prog-m1-d1-p02/cReg] {data/memory-saving-fmg.dat}; \addlegendentry{$p=2$}
        \addplot table[x=dofs-m1-d1-p03,y=prog-m1-d1-p03/cReg] {data/memory-saving-fmg.dat}; \addlegendentry{$p=3$}
        \addplot table[x=dofs-m1-d1-p04,y=prog-m1-d1-p04/cReg] {data/memory-saving-fmg.dat}; \addlegendentry{$p=4$}
        \addplot table[x=dofs-m1-d1-p05,y=prog-m1-d1-p05/cReg] {data/memory-saving-fmg.dat}; \addlegendentry{$p=5$}

        \addlegendimage{/pgfplots/refstyle={plot:m2d1p6}}
        \addlegendentry{$p=6$}
        \addlegendimage{/pgfplots/refstyle={plot:m2d1p7}}
        \addlegendentry{$p=7$}

        \nextgroupplot[
          title = {Poisson 2D},
        ]
        \addplot table[x=dofs-m1-d2-p01,y=prog-m1-d2-p01/cReg] {data/memory-saving-fmg.dat};
        \addplot table[x=dofs-m1-d2-p02,y=prog-m1-d2-p02/cReg] {data/memory-saving-fmg.dat};
        \addplot table[x=dofs-m1-d2-p03,y=prog-m1-d2-p03/cReg] {data/memory-saving-fmg.dat};
        \addplot table[x=dofs-m1-d2-p04,y=prog-m1-d2-p04/cReg] {data/memory-saving-fmg.dat};
        \addplot table[x=dofs-m1-d2-p05,y=prog-m1-d2-p05/cReg] {data/memory-saving-fmg.dat};

        \nextgroupplot[
          title = {Biharmonic 1D},
          extra description/.code = {
            \draw[-latex] (1.05,0.1) -- node[font=\scriptsize,rotate=90,anchor=north] {higher is better} (1.05,0.9);
          }
        ]
        \pgfplotsset{cycle list shift = 2}
        \addplot table[x=dofs-m2-d1-p03,y=prog-m2-d1-p03/cReg] {data/memory-saving-fmg.dat};
        \addplot table[x=dofs-m2-d1-p04,y=prog-m2-d1-p04/cReg] {data/memory-saving-fmg.dat};
        \addplot table[x=dofs-m2-d1-p05,y=prog-m2-d1-p05/cReg] {data/memory-saving-fmg.dat};
        \addplot table[x=dofs-m2-d1-p06,y=prog-m2-d1-p06/cReg] {data/memory-saving-fmg.dat}; \label{plot:m2d1p6}
        \addplot table[x=dofs-m2-d1-p07,y=prog-m2-d1-p07/cReg] {data/memory-saving-fmg.dat}; \label{plot:m2d1p7}

        \nextgroupplot[
          ylabel = {Solution vector\\[1.5ex]Storage saving vs.\\state-of-the-art},
        ]
        \addplot table[x=dofs-m1-d1-p01,y=prog-m1-d1-p01/cReg] {data/memory-saving-vector.dat};
        \addplot table[x=dofs-m1-d1-p02,y=prog-m1-d1-p02/cReg] {data/memory-saving-vector.dat};
        \addplot table[x=dofs-m1-d1-p03,y=prog-m1-d1-p03/cReg] {data/memory-saving-vector.dat};
        \addplot table[x=dofs-m1-d1-p04,y=prog-m1-d1-p04/cReg] {data/memory-saving-vector.dat};
        \addplot table[x=dofs-m1-d1-p05,y=prog-m1-d1-p05/cReg] {data/memory-saving-vector.dat};

        \nextgroupplot[
        ]
        \addplot table[x=dofs-m1-d2-p01,y=prog-m1-d2-p01/cReg] {data/memory-saving-vector.dat};
        \addplot table[x=dofs-m1-d2-p02,y=prog-m1-d2-p02/cReg] {data/memory-saving-vector.dat};
        \addplot table[x=dofs-m1-d2-p03,y=prog-m1-d2-p03/cReg] {data/memory-saving-vector.dat};
        \addplot table[x=dofs-m1-d2-p04,y=prog-m1-d2-p04/cReg] {data/memory-saving-vector.dat};
        \addplot table[x=dofs-m1-d2-p05,y=prog-m1-d2-p05/cReg] {data/memory-saving-vector.dat};

        \nextgroupplot[
          extra description/.code = {
            \draw[-latex] (1.05,0.1) -- node[font=\scriptsize,rotate=90,anchor=north] {higher is better} (1.05,0.9);
          }
        ]
        \pgfplotsset{cycle list shift = 2}
        \addplot table[x=dofs-m2-d1-p03,y=prog-m2-d1-p03/cReg] {data/memory-saving-vector.dat};
        \addplot table[x=dofs-m2-d1-p04,y=prog-m2-d1-p04/cReg] {data/memory-saving-vector.dat};
        \addplot table[x=dofs-m2-d1-p05,y=prog-m2-d1-p05/cReg] {data/memory-saving-vector.dat};
        \addplot table[x=dofs-m2-d1-p06,y=prog-m2-d1-p06/cReg] {data/memory-saving-vector.dat};
        \addplot table[x=dofs-m2-d1-p07,y=prog-m2-d1-p07/cReg] {data/memory-saving-vector.dat};
      \end{groupplot}
    \end{tikzpicture}%
    }%

    \subfloat{%
      \ref*{legend:savings}%
    }%

    \caption{%
      Estimated memory savings of compact multigrid compared to discretization-error accurate \gls*{fmg} in progressive \gls*{bfp} precision~\cite{kohl_multigrid_2024} (current state-of-the-art) for varying polynomial degrees.
      Here, we calculate with the experimentally determined constants from \cref{tab:experiments-bits}.
      A marker is placed every second refinement level.%
    }%
    \label{fig:memory-comparison}%
  \end{figure}

  The estimated memory savings are displayed in \cref{fig:memory-comparison}.
  For problem sizes where even the memory capacity of an entire supercomputer becomes prohibitive~\cite{gmeiner_quantitative_2016}, savings of around an order of magnitude are attainable.
  For BFPMG, the storage cost is primarily determined by the size of the precision increment, whereas for \gls*{cmg}, the (constant) bit-width of the finest grid is most important.
  Hence, the savings in the \gls*{fmg} solver increase with $m$, the precision increment of residual and correction.
  This shows in \cref{fig:memory-comparison} in the steeper slope of the \gls*{fmg} solver for the biharmonic equation compared to Poisson's equation.
  In contrast, no clear trend is visible w.r.t.\ the polynomial degree $p$, because with increasing $p$ both the solution vector's precision increment but also the number of bits on the finest grid increase (cf.\ \cref{tab:experiments-bits}).
  Comparing Poisson 1D and 2D, the dimension has a small effect on the storage savings.
  The \gls*{fmg} solver achieves a slightly smaller reduction in 2D than in 1D because at the same number of \glspl*{dof} the number of levels is smaller, and it is the number of levels that determines the precision.
  However, the higher the dimension, the quicker $n_\ell$ reduces toward coarser grids, reducing the penalty of having to store additional coarse-grid data in the compact representation.
  Therefore, if one compares 1D and 2D at the same level number, and thus at the same level of accuracy, the savings are larger in higher dimensions.

  \subsection{Computation}

  The time complexity of our scheme is $\mathcal O(n_L)$ arithmetic operations.
  As in previous work \cite{mccormick_algebraic_2021,tamstorf_discretization-error-accurate_2021,kohl_multigrid_2024}, the arithmetic must be performed in progressively higher precision.
  \Cref{alg:cfas-opt-recursive,alg:fmg-residual-opt-recursive} contain the bulk of the overall work.
  Just as with a standard multigrid V-cycle, \cref{alg:cfas-opt-recursive} comprises one product with the system matrix, a smoothing step, and some grid transfers on each level.
  The residual computation, on the other hand, is more expensive than in a standard method because it requires prolongating all sections of the compact solution to the finest grid.
  Moreover, we directly restrict the residual to the entire grid hierarchy, adding one restriction per level to the overall work (which would normally be part of the V-cycle).

  \subsection{Communication}

  $\mathcal O(n_L \log n_L)$ bits are transferred between memory and compute unit, matching the already known progressive precision \gls*{fmg}~\cite{mccormick_algebraic_2021,tamstorf_discretization-error-accurate_2021,kohl_multigrid_2024}.
  The dominating term stems from the temporaries in \cref{alg:cfas-opt-recursive,alg:fmg-residual-opt-recursive}.
  While the loop transformation from \cref{sec:memory-optimization} reduces the storage complexity, it does not change the amount of communication.
  With or without optimization, all entries of the temporaries are written/read to/from memory.
  On the finest grid, this amounts to $n_L$ elements with $\mathcal O(\log n_L)$ bits each.
  However, another benefit of the optimizations is that they increase spatial and temporal locality of memory accesses.
  Potentially, this makes these accesses less expensive, e.g., on CPUs if the optimized temporaries fit entirely into cache.

  \section{Numerical experiments}
  \label{sec:experiments}

  In this section, numerical results demonstrate that the proposed compact \gls*{fmg} solver achieves discretization-error accuracy with linear storage complexity and a linear number of operations.
  Our implementation builds on the GMP~\cite{granlund_gmp_2020} and MPFR~\cite{fousse_mpfr_2007} libraries for arbitrary precision arithmetic, and the G+SMO library~\cite{mantzaflaris_overview_2020} for the \gls*{fem} discretization.

  \subsection{Setup}

  We consider the homogeneous Dirichlet problem for Poisson's equation in one and two dimensions and the biharmonic equation in one dimension.
  Specifically, we define $\Omega = (0, 1)$ and $f(x) \in L^2\left( \Omega^d \right)$, and seek to approximate $u(x) \in C^{2m}$ so that
  \begin{equation}%
    \label{eq:pde-poisson}%
    \begin{aligned}
      - \Delta u(x) &= f(x) & \text{in $\Omega^d$}\,, \\
      u(x) &= 0 & \text{on $\partial(\Omega^d)$}\,,
    \end{aligned}%
  \end{equation}%
  \begin{equation}%
    \label{eq:pde-biharmonic}%
    \begin{aligned}
      u''''(x) &= f(x) & \text{in $\Omega$}\,, \\
      u(x) = u'(x) &= 0 & \text{on $\partial\Omega$}\,,
    \end{aligned}
  \end{equation}
  where $d \in \{1, 2\}$, $m = 1$ in \cref{eq:pde-poisson}, and $d = 1$, $m = 2$ in \cref{eq:pde-biharmonic}.

  To allow assessment of the error, we manufacture the analytical solutions
  $u(x) = x(1-x) \cos\left(\frac\pi 2 x\right)$
  for $(m, d) = (1, 1)$ (Poisson's equation \cref{eq:pde-poisson}, 1D),
  $u(x_0, x_1) = x_0 ( 1  - x_0 ) \cos( \frac{\pi}{2} x_0 ) x_1 ( 1 - x_1 ) \cos( \frac{\pi}{2} x_1 )$
  for $(m, d) = (1, 2)$ (Poisson's equation \cref{eq:pde-poisson}, 2D), and
  $u(x) = 1 - \cos(2 \pi x)$
  for $(m, d) = (2, 1)$ (biharmonic equation \cref{eq:pde-biharmonic}, 1D),
  giving rise to corresponding analytical \gls*{rhs} functions $f(x)$.
  The model problems are discretized using B-spline finite elements with uniformly refined knot vectors (corresponding to uniform structured meshes in more classical finite element methods).
  As a result, the storage cost of the mesh is insignificant, and stencils and \gls*{rhs} can be efficiently assembled on-the-fly.
  The polynomial degree $p$ varies from $p = 1, \dots, 5$ for Poisson's equation and $p = 3, \dots, 7$ for the biharmonic equation.
  As described in \cref{sec:precision}, we fix the number of compact V-cycles $\nu$ to $\nu = 1$ and only vary the number of \gls*{ir} iterations $N$ per level in the \gls*{fmg} solver.

  To verify our implementation, we compute a reference solution $u_L$ using a standard \gls*{fmg} solver in \gls*{fp} arithmetic with \num{200} ($(m,d)=(1,1)$), \num{100} ($(m,d)=(1,2)$), and \num{250} ($(m,d)=(2,1)$) bit mantissas and using 30 V(2, 1)-cycles per level to achieve discretization error.
  The compact \gls*{fmg} algorithm is considered successful if the error of its computed approximation $\tilde u_L$ is within twice the discretization error, i.e.,
  \begin{equation}
    \label{eq:experiment-successful}
    \|\tilde u_L - u\|_{H^m} \leq 2 \|u_L - u\|_{H^m}\,,
  \end{equation}
  and the convergence order is at most $0.05$ below the theoretically optimal rate, i.e.,
  \begin{equation}
    \label{eq:experiment-successful-order}
    -\log_2 \left( \frac{\|\tilde u_L - u\|_{H^m}}{\|\tilde u_{L-1} - u\|_{H^m}} \right) \geq (p - m + 1) - 0.05\,,
  \end{equation}
  for all levels $L \geq 4$, ignoring the preasymptotic region.

  We perform a simple experimental optimization to determine the remaining parameters $N$ (number of \gls*{ir}
  iterations per level) and $b_1$, $b_2$, $b_3$, $b_4$ (the precision constants in \cref{tab:precisions}) of the
  compact \gls*{fmg} algorithm.
  Initially, we fix $b_1 = b_3 = b_4 = 30$, $N = 12$ (which is more than sufficient to recover the convergence criteria
  \cref{eq:experiment-successful,eq:experiment-successful-order}) and increase $b_2$ (starting from $b_2$ = 1) until
  \cref{eq:experiment-successful,eq:experiment-successful-order} hold for $L_\mathrm{max} = 8$ ($d=1$) and $L_\mathrm{max} = 5$ ($d=2$).
  This is done for each problem and for each polynomial degree $p$ individually.
  Once the criteria are satisfied, we add an additional bit to $b_2$ (to add some stability when going to finer grids later) and ensure that the convergence criteria are still met.
  We then fix $b_2$ and proceed to tune $b_1$, $b_4$, $b_3$, and lastly $N$ (starting from $N = 1$) in the
  same way (also adding one additional bit and \gls*{ir} iteration respectively).
  It is certainly possible that the constants depend on each other, and that a different ordering of this
  search procedure would lead to different results.
  We chose this specific ordering hoping to find the smallest possible $b_2$ as that is likely to be the one
  constant that has the most impact on the storage savings.
  However, we do not claim that this is an optimal way to choose the parameters.
  In any case, having optimized the constants for these coarse problems, we fix all the precisions and run the experiments up to
  levels $L_\mathrm{max} = 22$ ($d = 1$) and $L_\mathrm{max} = 11$ ($d = 2$).

  \subsection{Results}

  \begin{table}[tbp]
    \tablesize
    \centering
    \caption{%
      Precision constants $b_1$, $b_2$, $b_3$, $b_4$ (cf.\ \cref{tab:precisions}) used in the experiments.
      $b_1 + 2 b_2$ is the number of bits (including signs) per \gls*{dof} consumed by all stored vectors combined ($\compact u$, $\std r$, $\compact y$) on the finest grid, independent of $h$.
      For example, with $p=3$, the compact 1D Poisson solver stores \num{15} bits per \gls*{dof} on the finest grid, compared to a standard solver using \num{96} bits in IEEE single precision (accurate for $L \leq 5$) or \num{192} bits in double precision (accurate for $L \leq 11$), cf.\ \cref{fig:hockey-sticks}.%
    }%
    \label{tab:experiments-bits}%
    \begin{tabular}{ll*{7}{>{\raggedleft\arraybackslash}c}}%
      \toprule
      \multirow{2}{*}{Constant (quantity)} & \multirow{2}{*}{Equation} &  &    &  & $p$ &   &    &    \\
      \cmidrule{3-9}
                                                     &               &  1 &  2 & 3 & 4 &  5 &  6 &  7 \\
      \midrule
                                                     & Poisson 1D    &  5 &  5 & 7 & 8 &  9 & -- & -- \\
      $b_1$ ($\compact u$)                           & Poisson 2D    &  4 &  5 & 5 & 7 &  9 & -- & -- \\
                                                     & Biharmonic 1D & -- & -- & 4 & 6 &  8 & 11 & 12 \\
      \midrule
                                                     & Poisson 1D    &  3 &  4 & 4 & 4 &  5 & -- & -- \\
      $b_2$ ($\std r$, $\compact y$)                 & Poisson 2D    &  4 &  4 & 4 & 5 &  6 & -- & -- \\
                                                     & Biharmonic 1D & -- & -- & 4 & 4 &  5 &  5 &  6 \\
      \midrule
                                                     & Poisson 1D    &  2 &  4 & 6 & 7 & 11 & -- & -- \\
      $b_3$ (matrices, $\std f$ in \pr{fmgResidual}) & Poisson 2D    &  2 &  3 & 4 & 7 & 15 & -- & -- \\
                                                     & Biharmonic 1D & -- & -- & 2 & 2 &  2 &  3 &  3 \\
      \midrule
                                                     & Poisson 1D    &  2 &  2 & 2 & 2 &  4 & -- & -- \\
      $b_4$ (matrices, $\std z$ in \pr{cFas})        & Poisson 2D    &  2 &  2 & 2 & 2 &  2 & -- & -- \\
                                                     & Biharmonic 1D & -- & -- & 3 & 2 &  2 &  2 &  2 \\
      \bottomrule
    \end{tabular}%
  \end{table}

  \begin{table}[tbp]
    \tablesize
    \centering
    \caption{%
      Number of \gls*{ir} iterations $N$ per \gls*{fmg} level (cf.\ \cref{alg:cfmg}) used in the experiments.%
    }%
    \label{tab:experiments-n-its}%
    \begin{tabular}{ll*{7}{>{\raggedleft\arraybackslash}c}}%
      \toprule
      \multirow{2}{*}{Equation} &&&& $p$ &&  &    \\
      \cmidrule{2-8}
                    &  1 &  2 & 3 & 4 &  5 &  6 &  7 \\
      \midrule
      Poisson 1D    &  4 &  3 & 4 & 5 &  9 & -- & -- \\
      Poisson 2D    &  3 &  2 & 4 & 7 &  9 & -- & -- \\
      Biharmonic 1D & -- & -- & 6 & 4 &  5 &  5 & 11 \\
      \bottomrule
    \end{tabular}%
  \end{table}

  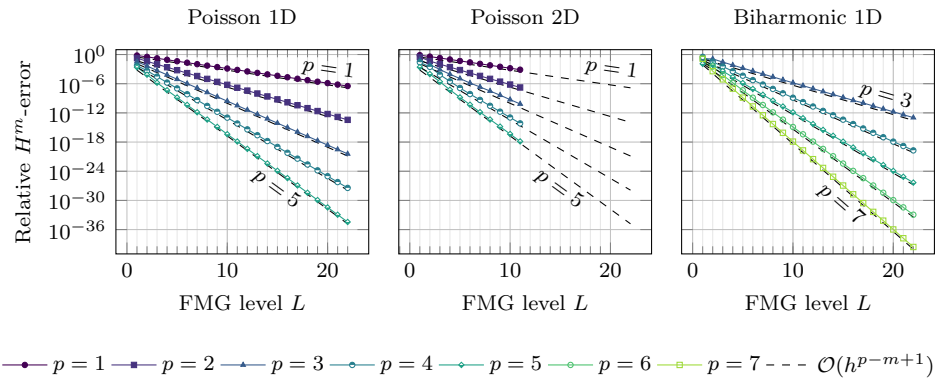
\begin{figure}[tbp]%
    \hfill%
    \subfloat{%
    \begin{tikzpicture}%
      \begin{groupplot}[
        group style = {group size=3 by 1, y descriptions at=edge left, horizontal sep=4mm},
        table/col sep = comma,
        width = 0.38\textwidth,
        every axis/.append style = {
          xlabel = {FMG level $L$},
          ymode  = log,
          ymax   = 10^1,
          ymin   = 10^-41,
          minor xtick = data,
          ytick  = {10^0, 10^-6, 10^-12, 10^-18, 10^-24, 10^-30, 10^-36},
        },
        table/x = level,
        cycle multiindex list = {[samples of colormap = 8 of viridis]\nextlist mymarklist},
        /tikz/label/.style = {sloped,node font=\footnotesize,black},
      ]
        \nextgroupplot[
          title = {Poisson 1D},
          ylabel = {Relative $H^m$-error},
          legend to name = {legend:results},
          legend columns = 8,
        ]
        \addplot table[x=level, y=H1] {data/poisson-lowprec-p01-d1-bfp.toml.csv}; \addlegendentry{$p=1$}
        \addplot table[x=level, y=H1] {data/poisson-lowprec-p02-d1-bfp.toml.csv}; \addlegendentry{$p=2$}
        \addplot table[x=level, y=H1] {data/poisson-lowprec-p03-d1-bfp.toml.csv}; \addlegendentry{$p=3$}
        \addplot table[x=level, y=H1] {data/poisson-lowprec-p04-d1-bfp.toml.csv}; \addlegendentry{$p=4$}
        \addplot table[x=level, y=H1] {data/poisson-lowprec-p05-d1-bfp.toml.csv}; \addlegendentry{$p=5$}

        \addlegendimage{/pgfplots/refstyle={pgf:p06}}
        \addlegendentry{$p=6$}
        \addlegendimage{/pgfplots/refstyle={pgf:p07}}
        \addlegendentry{$p=7$}

        \begin{scope}[domain=1:22,black,dashed,mark=none]
          \addplot {0.60*0.5^(1*x)} node [label,yshift=+1.7ex,pos=0.9] {$p=1$};
          \addplot {0.20*0.5^(2*x)};
          \addplot {0.10*0.5^(3*x)};
          \addplot {0.04*0.5^(4*x)};
          \addplot {0.02*0.5^(5*x)} node [label,yshift=-1.5ex,pos=0.7] {$p=5$};
        \end{scope}
        \addlegendentry{$\mathcal O(h^{p-m+1})$}

        \nextgroupplot[
          title = {Poisson 2D},
          minor xtick = {},
          minor x tick num = 9,
        ]
        \addplot table[x=level, y=H1] {data/poisson-lowprec-p01-d2-bfp.toml.csv};
        \addplot table[x=level, y=H1] {data/poisson-lowprec-p02-d2-bfp.toml.csv};
        \addplot table[x=level, y=H1] {data/poisson-lowprec-p03-d2-bfp.toml.csv};
        \addplot table[x=level, y=H1] {data/poisson-lowprec-p04-d2-bfp.toml.csv};
        \addplot table[x=level, y=H1] {data/poisson-lowprec-p05-d2-bfp.toml.csv};

        \begin{scope}[domain=1:22,black,dashed,mark=none]
          \addplot {0.60*0.5^(1*x)} node [label,yshift=+1.7ex,pos=0.9] {$p=1$};
          \addplot {0.20*0.5^(2*x)};
          \addplot {0.10*0.5^(3*x)};
          \addplot {0.04*0.5^(4*x)};
          \addplot {0.02*0.5^(5*x)} node [label,yshift=-1.5ex,pos=0.7] {$p=5$};
        \end{scope}

        \nextgroupplot[
          title = {Biharmonic 1D},
        ]
        \pgfplotsset{cycle list shift = 2}

        \addplot table[x=level, y=H2] {data/biharmonic-lowprec-p03-d1-bfp.toml.csv};
        \addplot table[x=level, y=H2] {data/biharmonic-lowprec-p04-d1-bfp.toml.csv};
        \addplot table[x=level, y=H2] {data/biharmonic-lowprec-p05-d1-bfp.toml.csv};
        \addplot table[x=level, y=H2] {data/biharmonic-lowprec-p06-d1-bfp.toml.csv}; \label{pgf:p06}
        \addplot table[x=level, y=H2] {data/biharmonic-lowprec-p07-d1-bfp.toml.csv}; \label{pgf:p07}

        \begin{scope}[domain=1:22,black,dashed,mark=none]
          \addplot {0.5*0.5^(2*x)} node [label,yshift=+1.7ex,pos=0.85] {$p=3$};
          \addplot {0.5*0.5^(3*x)};
          \addplot {0.5*0.5^(4*x)};
          \addplot {0.5*0.5^(5*x)};
          \addplot {0.5*0.5^(6*x)} node [label,yshift=-1.5ex,pos=0.7] {$p=7$};
        \end{scope}
      \end{groupplot}
    \end{tikzpicture}%
    }%
    \hfill\null%

    \hfill%
    \subfloat{%
      \ref*{legend:results}%
    }%
    \hfill\null%

    \caption{%
      Relative $H^m$-error $\|\tilde u_L - u\|_{H^m} / \|u\|_{H^m}$ of the computed compact solution $\tilde u_L$ to Poisson's equation \cref{eq:pde-poisson} ($m = 1$, $d \in \{1, 2\}$) and the biharmonic equation \cref{eq:pde-biharmonic} ($m = 2$, $d=1$) discretized with B-splines of degree $p$ (using parameters from \cref{tab:experiments-bits,tab:experiments-n-its}).
      The numerical results show the expected convergence rate $\mathcal O(h^{p - m + 1})$ consistent with theoretical results, where $h=2^{-L}$.
      It is worth comparing these results to \cref{fig:hockey-sticks} to get an idea of the limits of common IEEE floating-point formats.%
    }%
    \label{fig:experiment-convergence}%
  \end{figure}%

  The precision constants determined in this way are listed in \cref{tab:experiments-bits} for each equation and polynomial degree.
  Most are eight bits or lower, with exceptions occurring only for $p \geq 5$.
  The compact vectors may use as few as three bits for the finest section, one of which represents the sign.
  Notably, using linear elements for the 1D Poisson equation, as few as \num{5} bits for $\compact u_{L_\mathrm{max}}$ and \num{3} bits each for $\std r_{L_\mathrm{max}}$ and $\compact y_{L_\mathrm{max}}$ are required per \gls*{dof}.
  This means that for all vectors combined, only \num{11} bits per fine-grid point, about one sixth of a double-precision vector, must be stored to achieve discretization error for any $h$.
  While the values of $b_3$ and $b_4$ are also small, they are much less significant for the overall storage cost than $b_1$ and $b_2$.
  This is because $b_3$ and $b_4$ are the number of bits (including sign) used on the coarsest grid for stencil entries, $\std f$, and $\std z$.
  On the fine grids, the precision of these quantities is dominated by the level-dependent increment, and only few values are stored simultaneously.
  As expected, $b_1$ and $b_3$ depend on the equation and polynomial degree $p$
  because the lower discretization error inherent to higher-order methods must be balanced by lower quantization error.
  Conversely, we observe that $b_2$ and $b_4$ depend only mildly on $d$, $m$, and $p$.

  \Cref{tab:experiments-n-its} confirms that a low number of \gls*{ir} iterations is sufficient to achieve discretization-error accuracy.
  This indicates that compact \gls*{fmg} converges well despite low precision storage.
  Because the solver is configured to perform only one Gauss-Seidel sweep per level in each iteration, $N$ is also the total number of smoothing steps on the finest grid, each of which has essentially the same cost as a matrix-vector product.
  Note that worsening of the convergence rates with increasing spline degree~$p$ is expected~\cite{donatelli_symbol-based_2017}.

  Finally, to verify optimal $h$-convergence, \cref{fig:experiment-convergence} shows the relative $H^m$-error, $\|\tilde u_L - u\|_{H^m}/ \|u\|_{H^m}$ of the computed compact solution
  $\tilde u_L$ against the analytically known solution $u$.
  As can be seen, our algorithm achieves the asymptotically optimal convergence rate $\mathcal O(h^{p-m+1})$, which is
  limited by the discretization error.
  As expected, the convergence curves for Poisson's equation are essentially identical in 1D and 2D.
  We highlight that using seventh-order elements, we computed a discretization-error-accurate solution to the biharmonic equation with an $H^2$-error of less than \num{e-39} in as few as \num{12} bits for the solution vector on the finest grid.
  It is worth comparing these results to \cref{fig:hockey-sticks} to get an idea of the limits of common IEEE floating-point formats.

  \section{Conclusion}
  \label{sec:conclusion}

  In this work, we have reduced the space complexity of dis\-cret\-ization-error accurate multigrid from $\mathcal O(n \log n)$ to $\mathcal O(n)$.
  It is estimated that for state-of-the-art problems limited by memory capacity~\cite{gmeiner_quantitative_2016}, this corresponds to a memory footprint reduction of about an order of magnitude.
  The time complexity is the same as that of well-known multigrid methods in standard arithmetic.

  While a theoretical framework remains to be worked out, we have found that, in many cases, six bits or fewer provide sufficient precision to store the solution, residual, and correction on the finest grid.
  Contemporary hardware may not provide support for the exact bit widths we require, so in practice
  the data must be packed bit-wise into memory and converted to a supported data-type on load.
  Because the memory bandwidth of modern hardware lags behind the computational performance, the cost of this conversion is generally insignificant compared to the communication cost~\cite{grutzmacher_using_2023}.
  The coarsest grid, by contrast, may require bit-widths beyond the widest data-type implemented in hardware.
  The needed high-precision arithmetic can be efficiently emulated using low-precision integer arithmetic~\cite{ootomo_dgemm_2024}.

  A parallel implementation of our algorithms is out of scope of this article, but
  we presently do not see any reason why a parallel implementation with up to $\mathcal O(n / \log n)$ processes would not achieve linear memory complexity.
  Especially in the parallel setting, polyhedral compilers~\cite{feautrier_polyhedron_2011,bondhugula_practical_2008} could prove valuable for implementing the presented memory optimizations.

  \bibliographystyle{siamplain}
  \bibliography{cmg}

  \appendix

\end{document}


\maketitle

  \section{Review of multigrid}
  \label{sm:sec:multigrid}

  This section gives a brief introduction into existing multigrid theory.
  It familiarizes the reader with the multigrid variant that our compact algorithm is based on, which is \pgls*{fmg} method using \gls*{ir}~\cite{higham_accuracy_2002} and the \gls*{fas} to solve the residual equation on each level.
  For details, the reader is referred to~\cite{brandt_multigrid_2011,briggs_multigrid_2000}.

  \subsection{Multigrid V-cycle}

  The multigrid V-cycle is an iterative solver for elliptic \glspl*{pde} with $h$-independent convergence.
  In other words, the solution can be computed to a \emph{fixed} target accuracy with a number of V-cycles independent of the grid spacing.
  This sets it apart from stationary iterations and Krylov methods, whose performance typically degrades as $h$ decreases.
  There are two different formulations of multigrid: the \gls*{cs} and \gls*{fas}.
  For linear systems, \gls*{cs} and \gls*{fas} are equivalent.

  \subsubsection{Correction scheme}

  As outlined in \cref{sec:compact-representation}, we are solving an elliptic \gls*{pde} by discretizing it on a set of $L+1$ nested grids.
  Prolongation matrices $P_\ell \in \mathbb R^{n_\ell \times n_{\ell-1}}$ transfer vectors from level $\ell-1$ to level $\ell$.
  In this section, we write the discretized \gls*{pde} on refinement level $L$ as
  \begin{equation}
    \label{sm:eq:fine-grid-cs}
    A_L \std y_L = \std r_L\,.
  \end{equation}
  In the following, we refer to the exact solution of \cref{sm:eq:fine-grid-cs} as $A_L^{-1} \std r_L$ and an approximate computed iterate as $\std y_L$, avoiding the need of iteration indices.
  For example, the error is written as $\std e_L = A_L^{-1} \std r_L - \std y_L$.

  This error consists of smooth and oscillatory components.
  Key to the algorithm's efficiency is using two different methods to target these two errors individually.
  Both must be effective at reducing their respective part of the spectrum but typically perform poorly as a standalone solver.

  Reducing oscillatory error modes is the responsibility of the \emph{smoother}.
  Typical choices are relaxation methods of the form~\cite{falgout_generalizing_2004}
  \begin{equation}
    \label{sm:eq:smoother}
    \std y_L \gets \std y_L + M_L (\std r_L - A_L \std y_L)\,,
  \end{equation}
  where $M_L$ approximates $A_L^{-1}$.
  Examples include weighted Jacobi, Gauss-Seidel, and polynomial smoothers~\cite{baker_multigrid_2011}.

  After a few iterations of the smoother, the error is typically dominated by smooth components.
  Assuming that the remaining smooth error can be represented reasonably well on a coarser grid ($\std e_L \approx P_L \std d_{L-1}$), it can be well approximated by solving the residual equation $A_L \std e_L = \std r_L - A_L \std y_L$ restricted to the next coarser grid:
  \begin{equation}
    \label{sm:eq:coarse-grid-cs}
    \underbrace{R_L A_L P_L}_{A_{L-1}} \std d_{L-1} = R_L (\std r_L - A_L \std y_L)\,.
  \end{equation}
  Here, the restriction $R_\ell \in \mathbb R^{n_{\ell-1} \times n_\ell}$, $1 \leq \ell \leq L$, transfers vectors from one grid to the next coarser grid and is typically defined as $R_\ell := P_\ell^T$.

  Because \cref{sm:eq:coarse-grid-cs} is formulated on a coarser grid, it can typically be solved efficiently.
  To that end, note first that, it has essentially the same structure as the fine-grid problem \cref{sm:eq:fine-grid-cs} and, second, that smooth modes from the fine grid appear more oscillatory w.r.t.\ the coarse grid.
  Therefore, \cref{sm:eq:coarse-grid-cs} is best solved recursively, by applying a smoother and computing a correction from a yet coarser grid.
  The recursion is stopped once the system is small enough to be solved directly by a \emph{coarse-grid solver}.

  Once the correction $\std d_{L-1}$ is computed, it is prolongated and added to the solution:%
  \begin{equation}%
    \std y_L \gets \std y_L + P_L \std d_{L-1}\,.
  \end{equation}
  This reduces the smooth error components that contaminate $\std y_L$.
  It is common to apply additional post-smoothing steps after updating the solution.
  Generally, both the smoother and the coarse-grid solver are approximate.
  Moreover, the error left after smoothing is typically not perfectly in the range of $P$, meaning that the correction typically introduces new albeit fairly small additional oscillatory errors.
  Therefore, more V-cycles are usually required to attain satisfactory accuracy, with increasingly more as the grid becomes increasingly finer.

  \subsubsection{Full approximation scheme}

  \Gls*{fas} is most commonly used to solve nonlinear problems.
  Nevertheless, it has applications in linear problems as well.
  For example, one byproduct is an error estimator, which can be used for stopping criteria and adaptive grid refinement~\cite{brandt_multigrid_2011}.
  Moreover, it leads to $\tau$-extrapolation, which raises the approximation order without using higher-order elements~\cite{brandt_multigrid_2011}.
  For our purposes, the most important property of \gls*{fas} is that it gives accurate solutions to the fine grid equation, but on the coarse grids.

  \Gls*{fas} performs the same steps as \gls*{cs}, but with a different coarse-grid variable.
  Here, the coarse-grid variable is not the correction, but instead closely approximates the solution to the fine-grid equation.
  It is defined as
  \begin{equation}
    \label{sm:eq:coarse-grid-variable-fas}
    \std y_{L-1} := \tilde R_L \std y_L + \std d_{L-1}\,,
  \end{equation}
  where $\tilde R$ is a restriction operator defined on the primal space (unlike $R$ which is defined on the dual space), and $\std d_{L-1}$ is the same coarse correction as in \cref{sm:eq:coarse-grid-cs}.
  One typically has some freedom in choosing $\tilde R$ without sacrificing efficiency~\cite{brandt_multigrid_2011}.
  We comment on our choice in \cref{sec:compact-fas}.
  $\std y_{L-1}$ is found by (approximately) solving the FAS coarse-grid equation
  \begin{equation}
    \label{sm:eq:coarse-grid-fas}
    \begin{aligned}
      A_{L-1} \std y_{L-1} &= A_{L-1} \tilde R_L \std y_L + R_L (\std r_L - A_L \std y_L) \\
                           &= R_L \std r_L + \underbrace{(A_{L-1} \tilde R_L \std y_L - R_L A_L \std y_L)}_{\std \tau_L}.
    \end{aligned}
  \end{equation}
  This equation is found by replacing $\std d_{L-1}$ in \cref{sm:eq:coarse-grid-cs} with the expression given by \cref{sm:eq:coarse-grid-variable-fas} and rearranging.
  $\std \tau_L$ is known as the \emph{fine-to-coarse defect correction}, which modifies the coarse-grid equation so that at convergence its solution coincides with the (restricted) fine-grid solution $\tilde R_L A_L^{-1} \std r_L$.

  Having found a good approximation on the coarse grid (again through recursion), the fine-grid variable is updated according to
  \begin{equation}
    \label{sm:eq:fas-apply-correction}
    \std y_L \gets \std y_L + P_L \std d_{L-1} = \std y_L + P_L (\std y_{L-1} - \tilde R_L \std y_L)\,.
  \end{equation}

  \subsection{Full multigrid}

  As $h$ reduces, so does the discretization error.
  While V-cycles reduce the error by a fixed amount in linear time, achieving discretization error takes an increasing number of cycles.
  That is why either \gls*{cs} or \gls*{fas} is often best used in \pgls*{fmg} scheme, as we now explain.

  \Gls*{fmg} solves the \gls*{pde} on successively finer grids with an accuracy comparable to the respective discretization error.
  To avoid confusing the variables of the outer \gls*{fmg} solver and the inner V-cycle, we now consider the \gls*{pde} $\mathcal{L}u(x)=f(x)$.

  \subsubsection{Iterating on the original equation on each FMG level}

  In its classical form, as given in \cref{sm:alg:fmg}, the solution from the previous \gls*{fmg} level $P_L \std u_{L-1}$ is used as initial guess for the next finer grid.
  Because this initial guess is close to discretization error on the previous grid, the error contained within must only be reduced by a fixed factor to reduce it to discretization error on grid $L$.
  When measuring error in the energy norm, this factor is $2^{p-m+1}$, where $p$ is the polynomial degree of the discretization and $2m$ is the order of the \gls*{pde}~\cite{strang_analysis_2008}.
  Due to the $h$-independent convergence of multigrid, a fixed number of V-cycles achieves this error reduction.
  This makes \gls*{fmg} optimal in the sense that it is discretization-error accurate and its arithmetic cost is equivalent to just a few relaxation sweeps, meaning that it has linear time complexity in infinite precision.

  \begin{algorithm}[tbp]%
    \algosize%
    \caption{\Acrlong*{fmg} iterating on the original equation on each level.}%
    \label{sm:alg:fmg}%
    \paragraph{Input}
    \begin{tabular}{ll}
      $\mathcal L$       & continuous linear operator \\
      $f(x)$             & continuous \acrshort*{rhs} \\
      $L_\text{max}$     & index of finest grid \\
      $N$                & number of \acrshort*{ir}-iterations \\
      $\nu$              & number of V-cycles per \gls*{fmg} level \\
    \end{tabular}

    \paragraph{Output}
    Solution $\std u_{L_\text{max}}$.
    \unskip\medskip

    \begin{pseudo}[kw,fullwidth]%
      $(A_0, \std f_0) \gets \pr{discretize}(\mathcal L, f(x), L=0)$ \\
      $\std u_0 \gets A_0^{-1} \std f_0$ \\
      \kw{for} $L=1 \dts L_\text{max}$ \kw{step} $1$ \\+
        $(A_{0\rng L}, P_{1\rng L}, R_{1\rng L}, \std f_L) \gets \pr{discretize}(\mathcal L, f(x), L)$ \\
        $\std u_L \gets P_L \std u_{L-1}$ \\
        repeat $\nu$ times \\+
          $\std u_L \gets \pr{V-cycle}(A, P, R, \std u_L, \std f_L)$ \hfill\ct{\acrshort*{cs} or \acrshort*{fas}.} \\--
      \kw{return} $\std u_{L_\text{max}}$
    \end{pseudo}%
  \end{algorithm}

  \begin{algorithm}[tbp]%
    \algosize%
    \caption{\Acrlong*{fmg} using \acrlong*{ir} on each level.}%
    \label{sm:alg:fmg-correction}%
    \paragraph{Input}
    \begin{tabular}{ll}
      $\mathcal L$       & continuous linear operator \\
      $f(x)$             & continuous \acrshort*{rhs} \\
      $L_\text{max}$     & index of finest grid \\
      $N$                & number of \acrshort*{ir}-iterations \\
      $\nu$              & number of V-cycles per \gls*{fmg} level \\
    \end{tabular}

    \paragraph{Output}
    Solution $\std u_{L_\text{max}}$.
    \unskip\medskip

    \begin{pseudo}[kw,fullwidth]%
      $(A_0, \std f_0) \gets \pr{discretize}(\mathcal L, f(x), L=0)$ \\
      $\std u_0 \gets A_0^{-1} \std f_0$ \\
      \kw{for} $L=1 \dts L_\text{max}$ \kw{step} $1$ \\+
        $(A_{0\rng L}, P_{1\rng L}, R_{1\rng L}, \std f_L) \gets \pr{discretize}(\mathcal L, f(x), L)$ \\
        $\std u_L \gets P_L \std u_{L-1}$ \\
        repeat $N$ times \\+
          $\std r_L \gets \std f_L - A_L \std u_L$ \\
          $\std y_L \gets 0$ \\
          repeat $\nu$ times \\+
            $\std y_L \gets \pr{V-cycle}(A, P, R, \std y_L, \std r_L)$ \hfill\ct{\acrshort*{cs} or \acrshort*{fas}.} \\-
          $\std u_L \gets \std u_L + \std y_L$ \\--
      \kw{return} $\std u_{L_\text{max}}$
    \end{pseudo}%
  \end{algorithm}

  \subsubsection{Using IR on each FMG level}

  A different formulation of \gls*{fmg} is laid out in \cref{sm:alg:fmg-correction}.
  Instead of iterating on the original equation, it solves the residual/correction equation
  \begin{equation}
    \label{sm:eq:residual}
    A_L \std y_L = \std r_L := \std f_L - A_L P_L \std u_{L-1}
  \end{equation}
  on each \gls*{fmg} level, and then adds the correction to the solution.
  A long-established technique known as \gls*{ir} is to compute the residual in higer precision, solve for the correction, and update the solution in a loop with $N \geq 1$ iterations~\cite{higham_accuracy_2002,mccormick_algebraic_2021}.

  While equivalent to \cref{sm:alg:fmg} in infinite precision (if the total number of V-cycles is the same), this formulation facilitates reduced-precision arithmetic.
  Only the residual must be computed in relatively high precision, whereas the precision requirements on the V-cycle reduce significantly~\cite{mccormick_algebraic_2021,tamstorf_discretization-error-accurate_2021}.
  Concretely, in contrast to $\std u_L$ requiring progressive precision with $p+1$ bits increment, it is shown for $\nu=1$ that the V-cycle can be performed with just $m$ bits increment~\cite{tamstorf_discretization-error-accurate_2021}.

  \bibliographystyle{siamplain}
  \bibliography{cmg}